\documentclass[twocolumn]{autart}    

\usepackage{graphicx}          
\usepackage{eucal}
\usepackage{graphicx}
\usepackage{amsmath}
\usepackage{amssymb}
\usepackage{xcolor}
\usepackage{mdframed}
\usepackage{subfig}
\usepackage{lscape}

\begin{document}

\begin{frontmatter}

\title{Persistence of Excitation in  Reproducing Kernel Hilbert Spaces, Positive Limit Sets, and Smooth Manifolds} 


\author[VTME]{Andrew J. Kurdila}\ead{kurdila@vt.edu},    
\author[VTME]{Jia Guo}\ead{jguo18@vt.edu},               
\author[VTME]{Sai Tej Paruchuri}\ead{saitejp@vt.edu},    
\author[UMichAE]{Parag Bobade}\ead{paragsb@umich.edu}

\address[VTME]{Department of Mechanical Engineering, Virginia Tech, Blacksburg, VA 24060, USA}
\address[UMichAE]{Department of Aerospace Engineering, University of Michigan, Ann Arbor, MI 48109, USA}


\begin{abstract}                          
This paper studies the relationship between the positive limit sets of continuous semiflows and the newly introduced definition of persistently excited (PE) sets and associated subspaces of reproducing kernel Hilbert (RKH) spaces. It is shown that if the RKH space contains a rich collection of cut-off functions, persistently excited sets are contained as subsets of the positive limit set of the semiflow. The paper demonstrates how the new PE condition can be used to guarantee convergence of function estimates in the RKH space embedding method for adaptive estimation. In particular, the paper is applied to uncertain ODE systems with positive limit sets given by certain types of smooth manifolds, and it establishes convergence of adaptive function estimates over the manifolds.
\end{abstract}

\end{frontmatter}

\section{Introduction}

\vspace*{-.1in} 

In this paper we study the method of reproducing kernel Hilbert (RKH) space embedding for adaptive estimation of uncertain, or unknown, dynamic systems that are governed by systems of coupled, nonlinear ordinary differential equations (ODEs). The RKH embedding method for adaptive estimation has been introduced  in \cite{kl2013,bmpkf2017,bmpkf2017C}. This general formulation constructs estimates in $\mathbb{R}^d$ of the state of the unknown governing ODEs as well as estimates of an unknown function contained in the RKH space $ H$ that characterizes the uncertain governing ODEs.    
This paper investigates several unanswered questions related to the new notion of  persistency of excitation  (PE)  that has been introduced in the latter two of  these three papers. We derive  relationships  between the PE condition over an indexing set $\Omega$ that is a subset of the state space $X$ and the positive limit sets of semiflows over $X$. We also  construct or select good kernels that define the RKH space $H$ in applications where the governing semiflows exhibit certain asymptotic structural properties. These latter properties are expressed in terms of PE  conditions over  some classes of  smooth manifolds.  

The RKH embedding method generates a distributed parameter system, and its associated estimates  evolve in the generally  infinite dimensional space $\mathbb{R}^d \times H$.   There are many nontrivial questions about approximations and realizable implementations of the method, and  some study of the convergence of finite dimensional approximations of solutions of the RKH embedding equations is given  in \cite{bmpkf2017,bmpkf2017C}. In this short paper,  we only consider the convergence of the estimates generated by the governing  DPS system in the infinite dimensional state space $\mathbb{R}^d\times H$ that defines the RKH embedding formulation. The current investigation can be viewed as providing needed insight and intuition into the structure of the solutions of the RKH embedding equations, which is much needed for the  effective choice of approximating subspaces in practical implementations. 

\vspace*{-.1in} 
\subsection{Adaptive Estimation for Uncertain Nonlinear ODEs}

\vspace*{-.1in} 

A common setup for estimation of uncertain nonlinear systems  starts with an ordinary differential equation that can be decomposed into known and unknown parts, 
\vspace*{-.1in} 
\begin{equation}
\dot{x}(t)=g_0(x(t)) + g(x(t)), \quad x(0)=x_0 
\label{eq:orig_decomp}
\end{equation}

\vspace*{-.3in} 
with $x(t)\in \mathbb{R}^d$ for $t\in \mathbb{R}^+$,  $g_0:\mathbb{R}^d \rightarrow \mathbb{R}^d$ a known function,  and $g:\mathbb{R}^d \rightarrow \mathbb{R}^d$ an unknown function. One important problem  of adaptive estimation for such a  nonlinear system is to use the full state observations $\{x(t)\}_{t\in\mathbb{R}^+}$ to construct an evolution law for a state  estimate $\hat{x}(t)$ that approximates $x(t)$, in the sense that  $\hat{x}(t) \rightarrow x(t)$ as $t\rightarrow \infty$. In the language of adaptive estimation this is referred to as convergence of state  estimates.
A canonical model estimator for the original equation might choose the evolution law for the estimate to be 
\vspace*{-.1in} 
\begin{equation}
    \dot{\hat{x}}(t)=A\hat{x}(t)+g_0(x(t))+\hat{g}(t,x(t)) - Ax(t) \label{eq:est_decomp}
\end{equation}

\vspace*{-.3in} 
for a known matrix $A\in \mathbb{R}^d$, although many alternatives exist of course. Here $\hat{g}(t):=\hat{g}(t,\cdot)$ is an estimate of the unknown function $g$.  On defining  the state error $\tilde{x}(t)=x(t)-\hat{x}(t)$ and the function error $\tilde{g}(t):= \tilde{g}(t,\cdot):=g(\cdot)-\hat{g}(t,\cdot)$, the associated  error equation is obtained as 
\vspace*{-.1in} 
\begin{equation}
\dot{\tilde{x}}(t)=A\tilde{x}(t) + \tilde{g}(t,x(t)).
\end{equation}

\vspace*{-.3in} 
At a bare minimum then,   adaptive estimation methods for the above uncertain nonlinear ODEs must guarantee that trajectories of this error equation converge to zero.  
It is usually  considerably more difficult to guarantee that the    time-varying function estimate $\hat{g}(t,\cdot):\mathbb{R}^d \rightarrow \mathbb{R}$ converges in the sense that  $\hat{g}(t,\cdot)\rightarrow g$ as  $t\rightarrow \infty$.  It is this latter problem that is the primary  concern of this paper. 

\vspace*{-.1in} 
To gain some appreciation of the issues and nuances arising in the function estimation  problem, we consider two examples. Figures \ref{fig:circle} and  \ref{fig:fish} depict the phase portaits of the uncertain systems studied in Examples \ref{ex:circle} and  \ref{ex:fish}, respectively. These figures also include plots of the error in function estimates obtained by the RKH embedding techniques with the kernel of the RKH space selected as described in \ref{cor:cor2}. 
\vspace*{-.1in} 
\begin{exmp}
\label{ex:circle}
 The first example is a case of a supercritical Hopf bifurcation, which can be found in many textbooks on dynamical systems \cite{hale_kocak,khalil}. The system equations are given by
 \vspace*{-.1in} 
\begin{align}
\begin{Bmatrix} \dot{x}_1\\ \dot{x}_2
\end{Bmatrix}
&= \begin{Bmatrix} x_2 + x_1(1 - x_1^2 - x_2^2) \\
 -x_1 + x_2(1 - x_1^2 - x_2^2)
\end{Bmatrix} .
\label{eq:circle}
\end{align}

\vspace*{-.3in} 
Here we define $g_0(x):=\{x_2+x_1(1-x_1^2-x_2^2),0\}^T$, $g(x):=\{0,-x_1+x_2(1-x_1^2-x_2^2)\}$ in Equation \ref{eq:orig_decomp}. 
The figures below make clear that the positive limit set $\omega^+(x_0)$ is the circle $S^1$ for all $x_0\in \mathbb{R}^2$ for this dynamical system. When the method of RKH embedding is applied to this uncertain nonlinear system, we can obtain estimates $\hat{g}(t)$ of $g$ whose error is depicted in Figures \ref{fig:circle}(b,c).  These estimates have been constructed from finite dimensional approximations as discussed in \cite{bmpkf2017} using basis functions that are a collection of  (extrinsic) Sobolev-Matern kernels discussed in Corollary \ref{cor:cor2} centered over the positive limit set. 
\begin{figure*}[t]
    \centering
    \subfloat[Phase portrait]{\includegraphics[width=.3\textwidth]{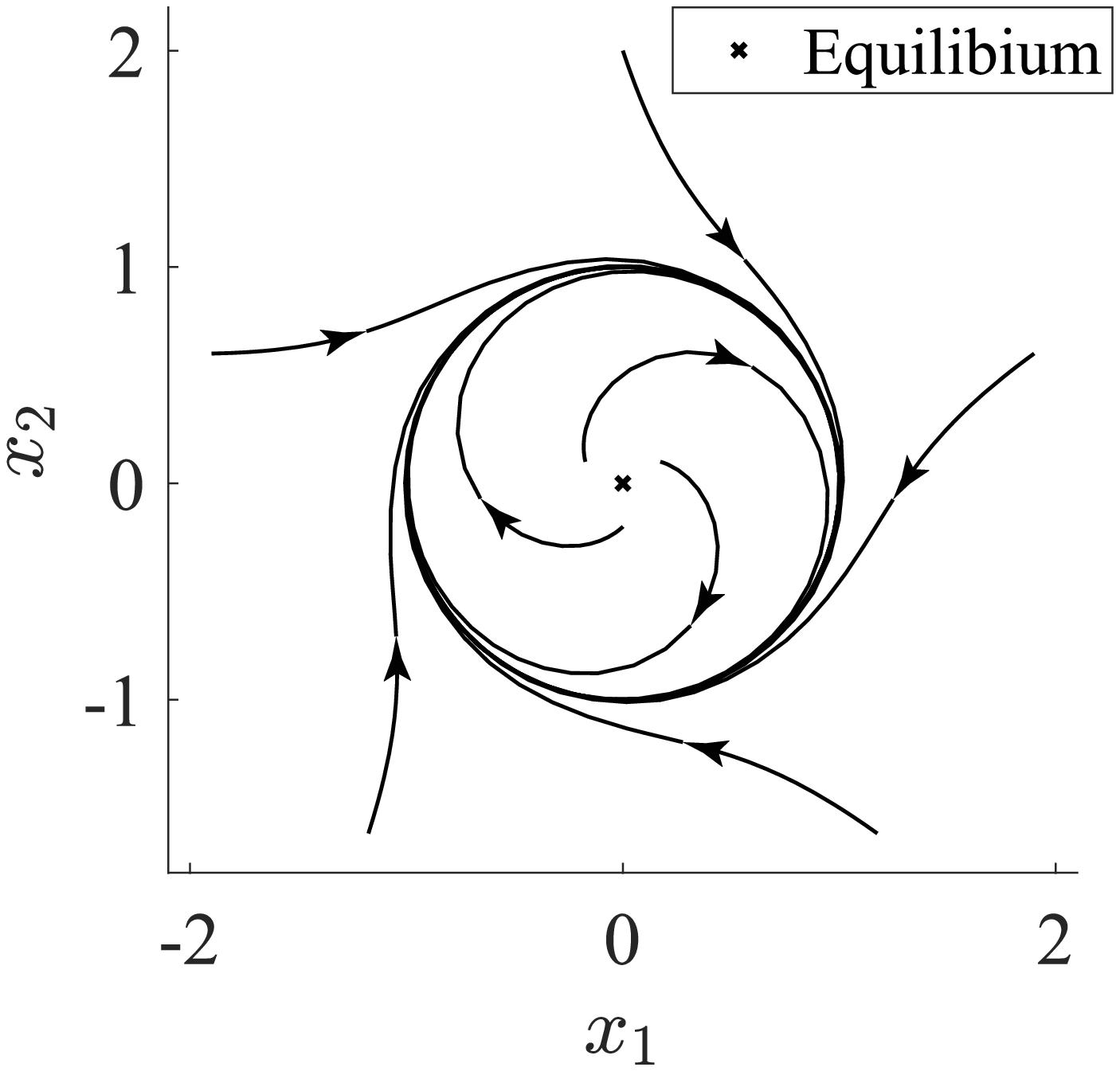}}
    \hspace{0.5cm}
    \subfloat[Error in Function Estimate]{\includegraphics[width=.3\textwidth]{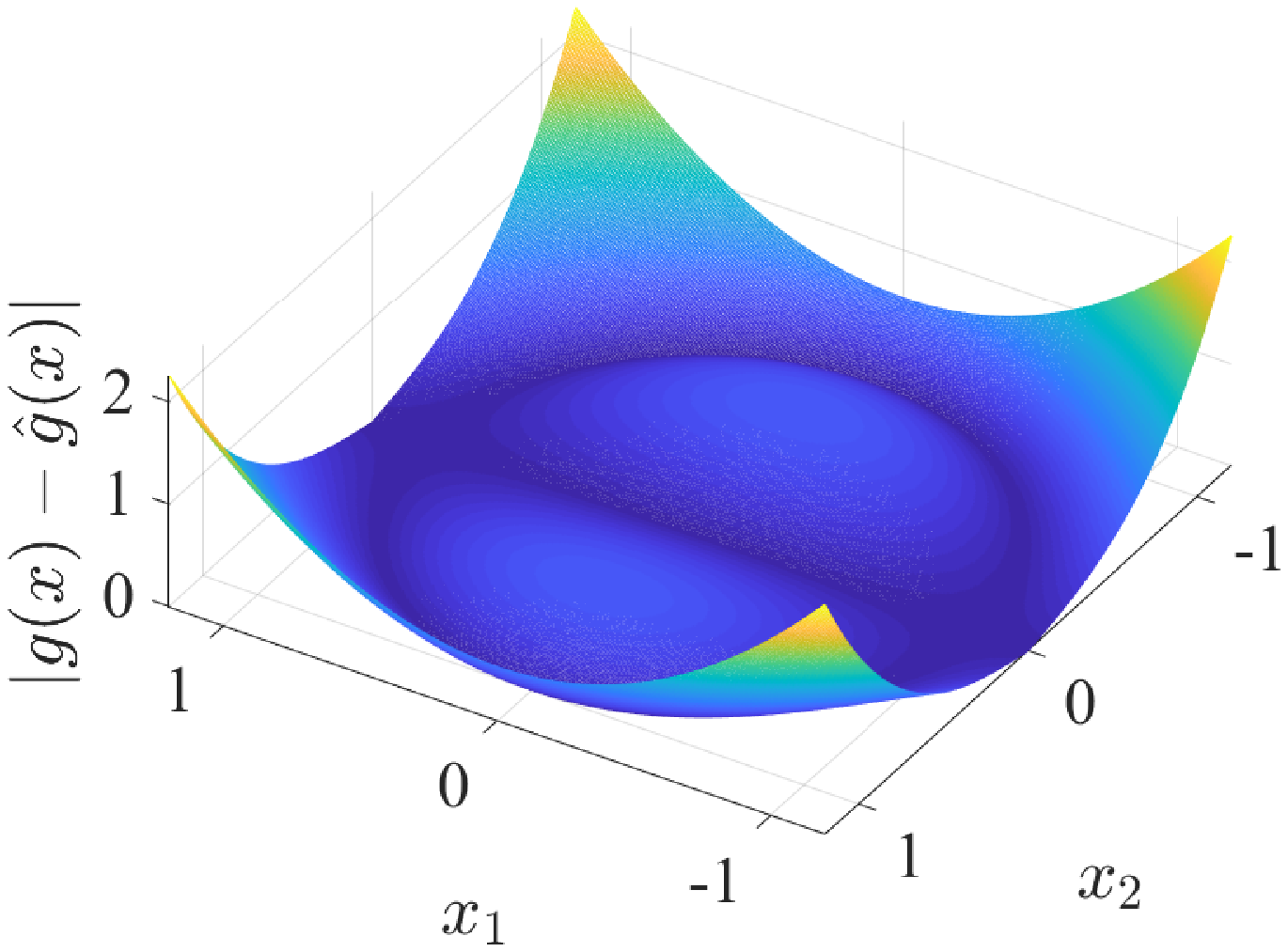}}
    \hspace{0.5cm}
    \subfloat[Error Contour]{\includegraphics[width=.3\textwidth]{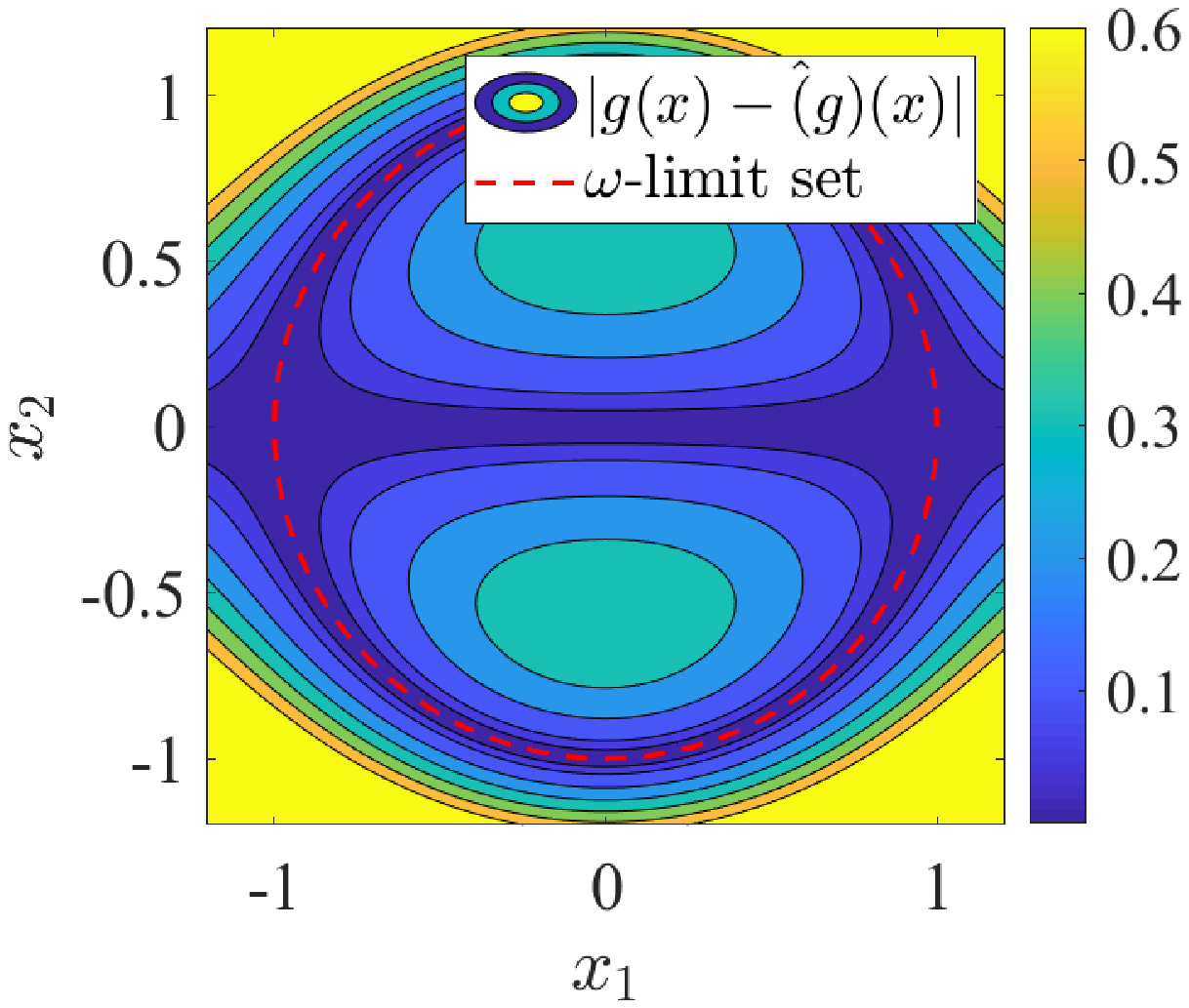}}
    \caption{Example 1}
    \label{fig:circle}
\end{figure*}


\end{exmp}
\vspace*{-.1in} 
\begin{exmp}
\label{ex:fish}
 In this example, the dynamical system contains a homoclinic loop. The example is studied in detail in \cite{hale_kocak}. The governing equations are 
 \vspace*{-.1in} 
\begin{align}
\label{eq:fish}
\begin{Bmatrix}
\dot{x}_1 \\ \dot{x}_2 \end{Bmatrix}
&= \begin{Bmatrix} 2x_2 \\ 
 2x_1 - 3x_1^2 + \lambda x_2(x_1^3 - x_1^2 + x_2^2)
 \end{Bmatrix}  
\end{align}

\vspace*{-.3in} 
In this example we define $g_0(x):=\{2x_2,0\}^T$ and $g(x):=\{0,2x_1-3x_1^2+\lambda x_2(x_1^3-x_1^2+x_2^2) \}^T$ in Equation \ref{eq:orig_decomp}. Again, application of the RKH embedding method of adaptive estimation to this problem can yield approximations $\hat{g}(t)$ of $g$ with error depicted in Figure \ref{fig:fish}(b,c). These estimates have been constructed from finite dimensional approximations as discussed in \cite{bmpkf2017} using basis functions that are a collection of  (extrinsic) Sobolev-Matern kernels discussed in Corollary \ref{cor:cor2} centered over the positive limit set.

\begin{figure*}[t]
    \centering
    \subfloat[Phase portrait]{\includegraphics[width=.3\textwidth]{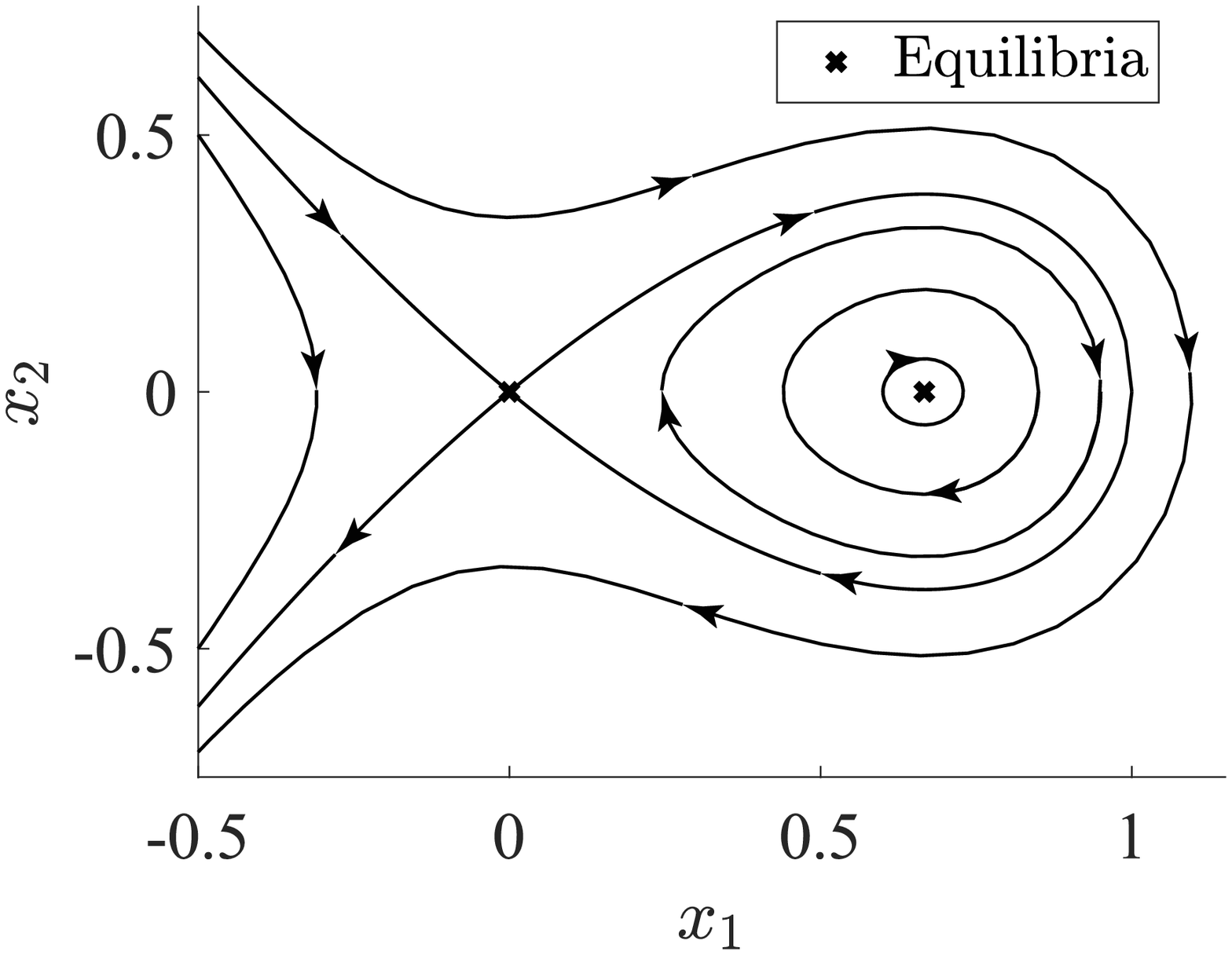}}
    \hspace{0.5cm}
    \subfloat[Error in Function Estimate]{\includegraphics[width=.3\textwidth]{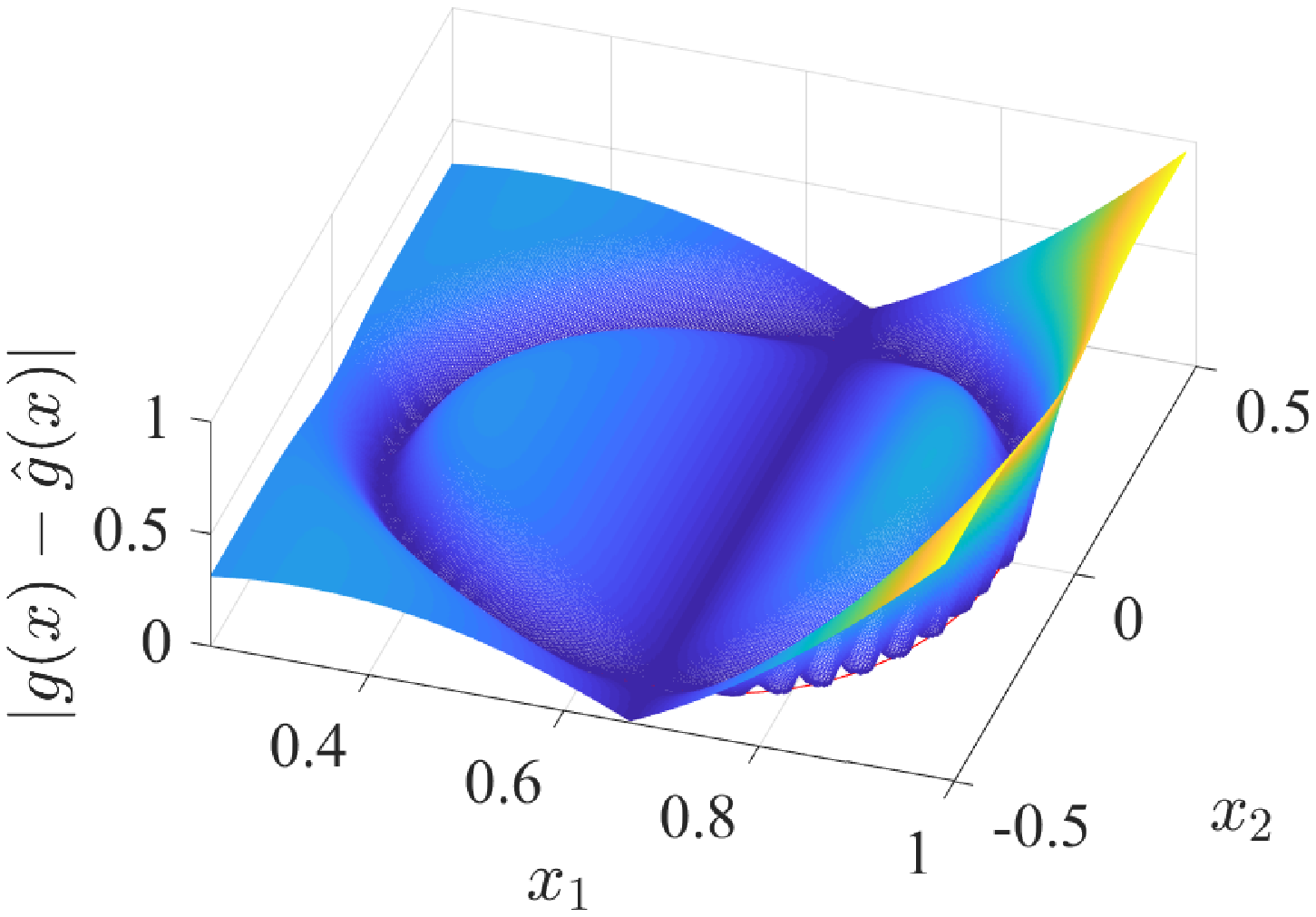}}
    \hspace{0.5cm}
    \subfloat[Error Contour]{\includegraphics[width=.3\textwidth]{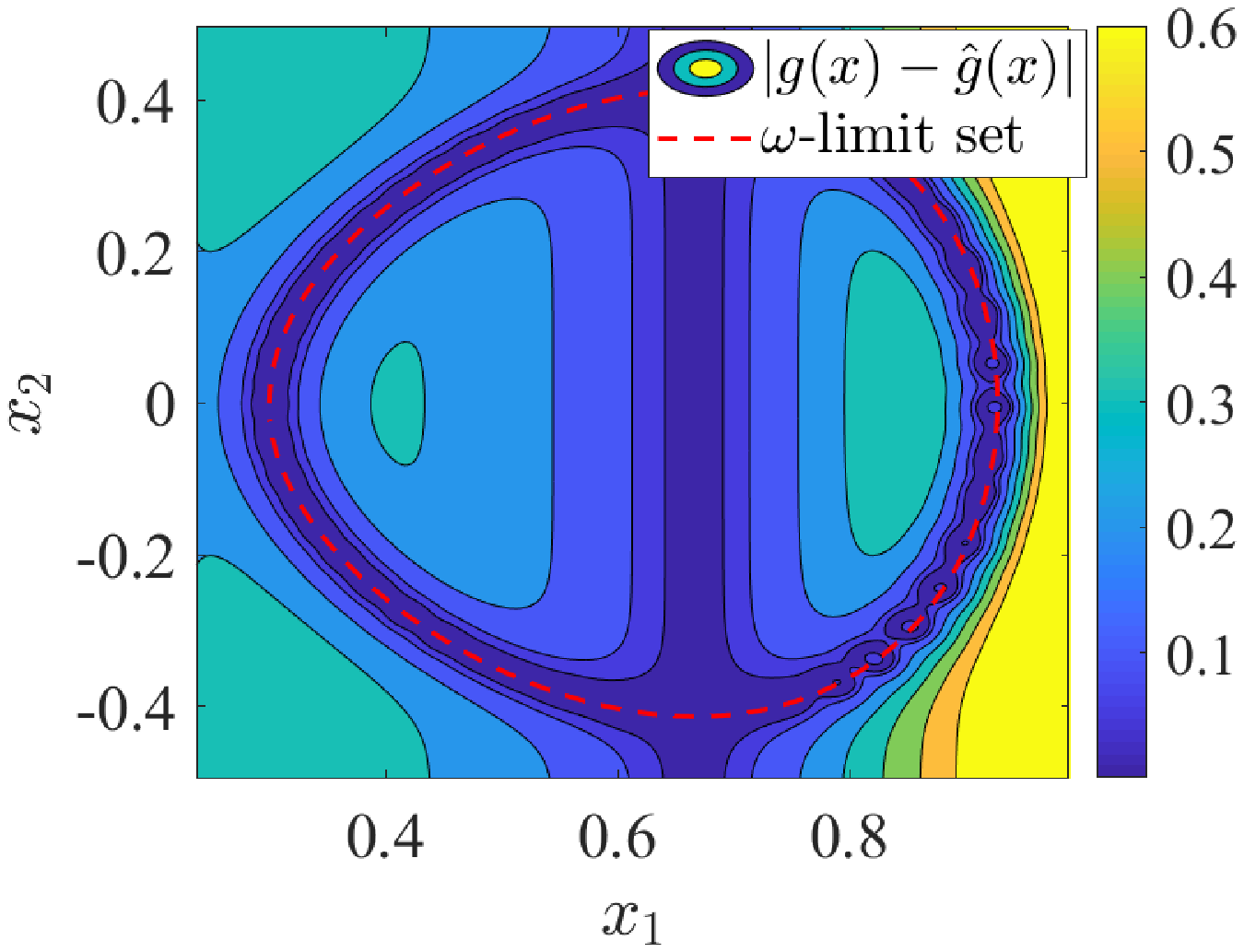}}
    \caption{Example 2}
    \label{fig:fish}
\end{figure*}
%
\end{exmp}

\vspace*{-.1in} 
Several observations about these two examples of application of the RKH embedding method are noteworthy and motivate this paper. 
In each case, the governing equations have the form of the nonlinear ODEs given in Equations \ref{eq:orig_decomp},  and are have error equations of the form in Equation \ref{eq:ode_Rd} that is studied in  detail in this paper. The first important observation to make from  the examples   is to note that the positive orbit $\Gamma^+(x_0):=\cup_{t\in\mathbb{R}^+}x(t)$ is the only data that is used to construct estimates of the unknown function. If we were interested in some offline, optimization-based estimate of an unknown function, it would come as no surprise that its estimates consist of functions that are supported on or near the prescribed data. We will see that, roughly speaking,  convergence of the RKH embedding method is guaranteed by a newly introduced  PE condition and  estimates of the unknown function  are built over regions of the state space where trajectories are in some sense ``concentrated.'' Here the notion of concentration is understood in terms of the positive limit set $\omega^+(x_0)$, which is known to attract trajectories if the orbit is precompact. \cite{walker}  

Moreover, both of   the positive limit sets in the examples  are striking and exhibit considerable structure: they can often be interpreted as manifolds. In Example \ref{ex:circle} shown in  Figure \ref{fig:circle}, global  solutions of the governing equations exist for every initial condition in $\mathbb{R}^d$. All trajectories converge to the positive limit set, which happens to be the canonical connected, compact, Riemannian manifold,  $S^1$. The flows generated for various parameters in $\lambda$ in Example \ref{ex:fish} exhibit more diverse qualitative limiting behaviors. As shown in Figure \ref{fig:fish}(a) when $\lambda=0$, the homoclinic loop encircles a stable region. Trajectories inside this region are all limit cycles. For these initial conditions, the positive limit sets are smooth,  regularly embedded submanifolds of $\mathbb{R}^d$. The form of these embedded manifolds is not as simple as in Example \ref{ex:circle}, that is, they are not one of the  well-known, ``iconic'' manifolds. When $\lambda<0$, the equilibrium $(x_e,0)$ for $x_e>0$ becomes unstable. It can be shown that the homoclinic loop becomes the $\omega$-limit set of all the trajectories starting from this region \cite{hale_kocak}.  

In either case, the examples illustrate a phenomenon that is common to many uncertain estimation problems. While the observations $\Gamma^+(x_0)=\cup_{t\in \mathbb{R}^+}x(t)$ are contained in $\mathbb{R}^d$, there is an underlying set or manifold that supports, approximately supports, or attracts the observed trajectories. We are interested in this paper in understanding conditions that establish that the RKH embedding method ``converges over'' these underlying structures. 

To frame our discussion of the RKH embedding method we briefly review the general strategy of ``linear-in-parameters'' (LIP) methods for adaptive estimation of uncertain nonlinear systems of ODEs. So-called LIP estimation might best be described as a part of the technical folklore for methods in adaptive estimation. 
This approach is ubiquitous in the adaptive estimation literature and is a  well-known tool among researchers who study this topic.
It   is safe to say that the most popular versions of adaptive estimation for the above type of uncertain  nonlinear ODEs choose  the function estimate in terms of a  linear-in-parameters  representation 
$
\hat{g}(t,\cdot)=\sum_{k=1}^n
\phi_k(\cdot)\alpha_k(t)=\Phi^T(\cdot)  \alpha(t)
$
with $\alpha_k(t)$ a time-varying parameter,  $\phi_k(\cdot):\mathbb{R}^d\rightarrow \mathbb{R}^d$ a function for $k=1,\ldots,n$, the vector $\alpha(t):=\{\alpha_1(t),\ldots, \alpha_n(t)\}^T\in \mathbb{R}^n$, and the matrix of functions  $\Phi^T(\cdot)=[\phi_1(\cdot),\ldots,\phi_n(\cdot)]$. Here the functions in $\Phi$ are known as the regressors, a common term arising from applications in  nonlinear regression. 
If the unknown function $g$ has the representation $g=\sum_{k=1}^{n} \phi_k(\cdot)\alpha^*_k=\Phi^T(\cdot)\alpha^*$ for some unknown constants $\{\alpha^*_k\}_{k\leq n}$, then the error in function estimates is $\tilde{g}(t,\cdot)=\Phi(\cdot)\tilde{\alpha}(t)$ with the parameter error  $\tilde{\alpha}(t):=\alpha^*-\alpha(t)\in \mathbb{R}^n$.  In this case the error in the function estimates $\tilde{g}(t)\rightarrow 0$ if the finite set of parameters errors converge $\tilde{\alpha}(t)
\rightarrow 0\in \mathbb{R}^n$. It is for this reason that the task of estimating functions in the usual LIP framework reduces to questions of parameter convergence in $\mathbb{R}^n$.    

One of the foundations   of modern adaptive estimation for ODEs has been recognition of the fact that persistency of excitation conditions can be sufficient to guarantee parameter convergence. The notion of persistence of excitation in its conventional form, that is, as it pertains to the the ODE error Equations \ref{eq:ode_Rd}, is defined next.

\begin{defn}
\label{defn:PE_Rd}
The  regressors $\Phi$  are persistently excited by the positive orbit  $\Gamma^+(x_0)$ if there are positive constants $\gamma_1,\gamma_2,T,$ and $\Delta$ such that for each $t\geq T$,
\vspace*{-.1in} 
{\small 
\begin{equation}
\label{eq:PE_Rd}
    \gamma_1 \| \alpha \|^2_{\mathbb{R}^d} \leq 
    \int_t^{t+\Delta} \left( \alpha^{\top} \Phi\bigl( x(\tau) \bigr) \Phi\bigl( x(\tau)                \bigr)^{\top}\alpha \right) \; d\tau \leq
    \gamma_2\|\alpha\|^2_{\mathbb{R}^d}
\end{equation}
}

\vspace*{-.3in} 
for all $\alpha \in \mathbb{R}^N$.
\end{defn}
 The papers \cite{narkud,shimkin1987Persistency,naranna87,Moore1992Functional,Boyd1983On} and a  number of standard texts \cite{sb2012,naranna,IaSu,PoFar} on adaptive estimation make a careful study of this condition and how it facilitates a proof that the parameter error $\tilde{\alpha}(t)$ converges to zero as $t\rightarrow \infty$. 
In some cases it is too much to hope that all the parameters $\hat{\alpha}(t):=\{\hat{\alpha}_1(t), \ldots,\hat{\alpha}_n(t)\}$ in the approximations $\hat{g}(t,\cdot)=\sum_{k=1,\ldots,n}\phi_k(\cdot)\hat{\alpha}_k(t)$ converge. A  means of weakening the above PE condition introduces the notion of {\em partial} persistency of excitation.  
One version of the definition of a partial PE condition modifies  the inequalities above and  replaces them with the condition that
\vspace*{-.325in} 
{\small 
$$
\gamma_1\|P_V v\|_{\mathbb{R}^n}^2
\leq  \alpha^T \int_{t}^{t+\Delta} \Phi(x(\tau))\Phi^T(x(\tau)) d\tau \cdot \alpha 
\leq \gamma_2 \|P_Vv\|_{\mathbb{R}^n}^2
$$
}

\vspace*{-.3in} 
\noindent with  $P_V:\mathbb{R}^n \rightarrow \mathbb{R}^n$ a projection onto a linear subspace $V\subset \mathbb{R}^n$. This generalization then can be used to guarantee, as a special case, that only certain of the coefficient estimates converge, but not all. 

\vspace*{-.1in} 
As we will discuss in more detail shortly, the method of RKH embedding recasts the above adaptive estimation problem so that the state errors $\tilde{x}(t)$  and function errors $\tilde{g}
(t):=\tilde{g}(t,\cdot)$ evolve in a product space having the form $\mathbb{R}^d \times \mathbb{H}$ with $\mathbb{H}=H^d$ a vector-valued RKH space of functions. The space $H$ is known as the hypothesis space and its selection is based on what class of priors or information seems relevant regarding the estimation problem at hand.  
 The precise form of the PE definition in this paper, and the associated theorems that depend on it, are written for a model problem with the vector-valued function $g:=Bf$
with $f:\mathbb{R}^d \rightarrow  \mathbb{R}$ a scalar valued function and $B\in \mathbb{R}^{d\times 1}$. This restriction does not seem too severe, simplifies the notation considerably, and conveys the underlying geometric relationships between orbits $\Gamma^+(x_0)$ of semiflows starting at $x_0\in X$, persistency of sets $\Omega$,  and RKH spaces $H_\Omega$. Moreover, the  extension to general vector-valued functions would proceed in principle along the same lines as the strategy in \cite{bobade2} used for consensus estimation. 

\vspace*{-.1in} 
\subsection{Overview of New Results}

\vspace*{-.1in} 

In either of the papers \cite{bmpkf2017,bmpkf2017C} some of the standard questions regarding RKH embedding have been discussed such as existence of solutions and well-posedness, continuous dependence on initial conditions, as well as stability and convergence of finite dimensional approximations.  In this paper we focus primarily on building more intuition and insight regarding the newly introduced notion of persistency of excitation in the RKH embedding method. Starting with an RKH space $H_X=\overline{\text{span}\{\mathfrak{K}(x,\cdot)| x\in X\}}$ of functions over $X$, we then define for some {\em indexing set} $\Omega \subseteq X$ the subspace 
$H_\Omega:=
\overline{\text{span}\{ \mathfrak{K}(x,\cdot) \ | \ x\in \Omega \}}
$. Note carefully that functions in $H_\Omega$ are supported on $X$, which is why $\Omega$ is referred to as the indexing set. This space is related to, but distinct from, the space $R_\Omega(H_X)$ that are {\em restrictions} of functions in $H_X$ to the subset $\Omega$. 
We have the following definition of persistence of excitation for the RKH error Equations  \ref{eq:rkhs_err}.
\vspace*{-.1in} 
\begin{defn}
\label{defn:PE_rkhs}
The indexing set $\Omega$ and RKH space $H_\Omega$  are persistently excited by the orbit $\Gamma^+(x_0)$ if there are positive constants $\gamma_1,\gamma_2,T$, and $\Delta$ such that for each $t\geq T$,
\vspace*{-.2in} 
\begin{equation}
\label{eq:PE_RKHS}
    \gamma_1\|f\|^2_{H_\Omega} \leq
    \int_t^{t+\Delta} \left( \mathcal{E}_{x(\tau)}^*\mathcal{E}_{x(\tau)}f,f                    \right)_{H_X} \; d\tau \leq
    \gamma_2\|f\|_{H_\Omega}^2
\end{equation}
\end{defn}

\vspace*{-.1in} 
\noindent
Here $\mathcal{E}_x:f\mapsto f(x)$ is the evaluation functional at $x$ and $\mathcal{E}_x^*$ is its adjoint operator. 
The classical definition given above defines persistency of excitation for a specific set of regressors and the trajectory of a semidynamical system. The new PE condition holds for an indexing  set $\Omega\subseteq X$, space of functions $H_\Omega$, and a trajectory of a semidynamical system.
  It should be noted that the  PE condition above is over a set $\Omega\subset X$, which may or may not be the entire state space $X$. It this sense it bears some resemblance of the  partial  PE conditions that are defined over subspaces of parameters in $\mathbb{R}^n$. The similarity in form is all the more apparent when we note that $\|\cdot \|_{H_\Omega}:=\|P_\Omega(\cdot)\|_{H_X}$ with $P_\Omega$ the $H_X$-orthogonal projection onto $H_\Omega$: the closed subspace $H_\Omega$ is endowed with the norm it inherits  from $H_X$. 
It should also be pointed out that both the set $\Omega$ and the kernel $\mathfrak{K}_X$ (that determines $H_X$, and therefore determines $H_\Omega\subseteq H_X$)  are free to be selected when trying to apply the  above PE condition in the method of RKH embedding. 

\vspace*{-.1in} 
Intuitively, we expect some kernels are more useful than others in the RKH embedding method, and one of the primary thrusts of this paper is to explore the alternatives.  As we will see, we obtain a strong conclusion about what type of indexing sets $\Omega$ are PE when we restrict attention to kernels that define function spaces that are good at separating important subsets of $X$. There are many ways to think about how well the functions in a space $H_X$ separates points or sets. We find that  one important class of RKH spaces consists of functions that feature a rich set of (possibly smooth) cut-off or bump functions. Our first primary result in Theorem \ref{th:th2} is that if the RKH space does indeed contain a rich collection of these functions, then we have the following  implication, 
\begin{equation}
``\Omega \text{ and }  H_\Omega 
\text{ are PE''} \quad \Longrightarrow \quad \Omega \subseteq \omega^+(x_0),
\label{eq:PE_implies}
\end{equation}
with $\omega^+(x_0)$ the positive limit set of a trajectory starting at $x_0$. 
In other words in terms of the new definition of PE,  if a trajectory $\Gamma^+(x_0)$ persistently excites a  set $\Omega$, the set $\Omega$  is contained in the positive limit set $\omega^+(x_0)$.
This result provides novel insight into the structure of this type of persistently excited systems : PE sets are not transient but rather consist of points whose neighborhoods are visited by the trajectory  infinitely often. In fact, a bit more is actually required as illustrated in Theorem \ref{th:pe_pt}: the ``time of visitation'' is bounded below in a certain sense. This intuition should be compared with the interpretations of the usual Definition \ref{defn:PE_Rd}: a vector signal $t\mapsto \Phi(x(t))\in \mathbb{R}^n$ is (partially) PE if on average it visits all directions in (a subspace of) $\mathbb{R}^n$. 

\vspace*{-.1in} 
We should also emphasize at this point that while the intent of this paper is to inform and enhance our understanding of the RKH embedding method, the result in Equation \ref{eq:PE_implies} is not dependent on the fact that the trajectory  under study $t\mapsto x(t)\in \mathbb{R}^d$ happens to be the solution of our model ODE problem in Equation \ref{eq:est_decomp}. We have worked to express the condition in Equation \ref{eq:PE_implies} in very general terms. As long as the PE condition holds under the hypothesis described above, and  $\Gamma^+(x_0)$ is the forward orbit  of a continuous semiflow on the complete metric space $(X,d_X)$, we conclude that $\Omega \subseteq \omega^+(\Omega)$.

\vspace*{-.1in} 
It is  then natural to ask how to choose kernels that exhibit the separation properties that enable the conclusion above.
One approach, which we refer to as an intrinsic method, applies to cases in which  $X=\Omega$ is in fact a compact, connected, smooth, Riemannian manifold $M$.  Here we assume that $M$ is known and that the kernel over $M$ is known. Example \ref{ex:circle} is the type of problem we have in mind here, where the positive limit set is a simple well-known manifold. Numerous intrinsic kernels can be defined over the circle, the sphere, or more generally homogeneous manifolds \cite{HNW}.  In this case we choose kernels that guarantee that the native space $H_M$ is in fact equivalent to a certain Sobolev space $W^{r,2}(M)$ for $r$ large enough. That such equivalences are possible follows from the Sobolev embedding theorem. \cite{af2003} The Sobolev spaces defined over such a manifold $M$ can be shown to contain a rich family of smooth cutoff functions. In this framework, if the forward orbit $\Gamma^+(x_0)$ for some $t_0\in \mathbb{R}$ of  any continuous flow on $M$ is PE, Corollary \ref{cor:cor1} implies that the  manifold is transitive. That is, it supports a flow  that has a dense orbit. The study of when a particular manifold is transitive is of interest in its own right \cite{lopez}, so the new PE condition can be used to study whether a manifold is transitive. 

\vspace*{-.1in} 
While  this is an interesting result, it is not usually strictly or directly  applicable to understanding the convergence properties of the RKH embedding problem. There are two essential problems here. First, there are many problems where the positive limit set might be a nice smooth, compact, Riemannian  manifold $\omega^+(x_0)=M$, but we do not know the form of the manifold {\em a priori}. In such cases defining the kernel in closed form to be used in analysis or approximation is impossible. Example \ref{ex:fish} is of this type: the positive limit set is a smooth manifold, but it is not one over which catalogs of intrinsic kernels are defined. It is also possible,  on the other hand, that we do know the exact form of the manifold, but it is not one of the standard manifolds like the circle, sphere, or torus.  Even if  we in principle can define the kernel  through the fundamental solutions of certain elliptic differential operator on the manifold $M$ as in Corollary \ref{cor:cor1},  it may be intractable to compute this fundamental solution for the manifold at hand. This problem can be as hard, or harder perhaps, that the original estimation problem. 

\vspace*{-.1in} 
It should be kept in mind that the aim of the RKH embedding method is to carry out adaptive estimation of {\em uncertain} nonlinear ODEs. It is typically the case in such situations that the exact form of the positive limit set is unknown. That is, we are more interested in problems like Example \ref{ex:fish}, in contrast to Example \ref{ex:circle}.
In this case,  we assume the $M$  is an unknown, connected, smooth, 
 (regularly) embedded submanifold of $\mathbb{R}^d$, and  we resort to an extrinsic method. In this technique we build a well-defined kernel on $X$ for a large set $X$ that contains $M$, and then we define a kernel by restriction on the manifold $M \subseteq X$. It can be the case that a  plethora of kernels exist for good kernels over the large space $X=\mathbb{R}^d$.  Taking care to choose  the kernel smooth enough, we obtain a kernel on $M$ defined by restriction. The expression for the kernel on $M$ is given in terms of the kernel on the larger space $X$, which is known. Corollary \ref{cor:cor2} then shows that that $M=\omega^+(x_0)$ in this case. All of the numerical examples depicted in Figures \ref{fig:circle} and \ref{fig:fish} have been computed using this extrinsic method.

\vspace*{-.125in}
\section{Notation}
\vspace*{-.1in} 
In this paper the symbols $ \mathbb{N}^+, \mathbb{R},\mathbb{R}^+$ denote the non-negative integers, real numbers, and non-negative real numbers, respectively. The expression 
$
a \lesssim b 
$ 
means that there is a  constant $c>0$ that does not depend on $a,b$ such that $  a \leq c \cdot b $. The symbol $\gtrsim$ is defined similarly.  
The paper makes use of Lebesgue spaces and Sobolev spaces on subsets $\Omega$ of $\mathbb{R}^d$,  and it also uses these spaces when they are defined more generally on  measurable subsets $\Omega$ of  certain  Riemannian  manifolds $M$.   The norm on the Banach spaces
  $L^p(\Omega):=L^p_\mu(\Omega)$ of $\mu$-integrable  functions over  $\Omega \subseteq \mathbb{R}^d$ take the familiar  form 
$
\|f\|^p_{L^p(\Omega)}:=\int_\Omega |f(x)|^pd\mu
$
with the measure $\mu$ on $\mathbb{R}^d$ for  $1\leq p<\infty$, with the usual modification for $p=\infty$. Recall that  when $\Omega \subseteq \mathbb{R}^d$, the Sobolev space $W^{r,p}(\Omega)$ for a positive integer $r$ consists of functions that have weak derivatives of all orders less than or equal to $r$ in $L^p(\Omega)$, and the norm on these Banach spaces is usually written 
\vspace*{-.125in} 
\begin{equation}
\|f\|_{W^{r,p}(\Omega)}^p:=\sum_{0\leq |\alpha| \leq r} \int_\Omega \left |\frac{\partial^{|\alpha|}f}{\partial x^\alpha} \right |^p dx
\label{eq:sob_Rd}
\end{equation}

\vspace*{-.2in} 
with the summation taken over all multi-indices $\alpha=(\alpha_1,\ldots,\alpha_d)$,  $|\alpha|=\sum_{i=1,\ldots,d}\alpha_i$, and here the measure $\mu$ is selected to be Lebesgue measure $dx$.  The Sobolev spaces for non-integer $r>0$ are defined in terms of interpolation theory as discussed in \cite{af2003}. 
A bit more detail is required to define the spaces $L^p(\Omega)$ and $W^{r,p}(\Omega)$ for $\Omega\subseteq M$, with $M$ a manifold.  In this paper $M$ is always assumed to be a connected, complete Riemannian manifold with a positive injectivity radius and bounded geometry.  See \cite{triebel2}, Chapter 7 or \cite{HNW,HNRW} for a discussion of these properties.  For purposes in this paper, it suffices to note that  compact, connected Riemannian manifolds and $\mathbb{R}^d$ satisfy these conditions. For such  a Riemannian manifold $M$ denote the  metric $g$ and inner product $<\cdot,\cdot>_{g,p}$ on the tangent space $T_pM$. We define the associated volume measure $d\mu$ on $M$, and its local representation in terms of the set of coordinates $(x^1,\ldots,x^d)$ is given by $d\mu(x):=\sqrt{det(g)}dx^1\ldots dx^d$. The norm $\|f\|_{L^p(\Omega)}$ has the same expression given above with the measure selected to be the usual volume measure on the manifold $M$. 
The Banach spaces $W^{r,p}(\Omega)$ for measurable subsets $\Omega\subseteq M$ are equipped with the norm 
\vspace*{-.125in} 
\begin{equation}
\|f\|^p_{W^{r,p}(\Omega)}:=
\sum_{j=0,\ldots,r}\int_\Omega
|\nabla^jf|^p_{g,p}d\mu(p)
\label{eq:sob_man}
\end{equation}

\vspace*{-.2in} 
for $1\leq p<\infty$ where $\nabla$ is the covariant derivative over $(M,g)$. When applied to a set $\Omega\subseteq M=\mathbb{R}^d$, the definitions over the manifold $M$ define norms that are equivalent to the usual ones for Sobolev spaces defined on subsets of $\mathbb{R}^d$. As discussed in \cite{HNW,HNRW} in this case the expression in Equation \ref{eq:sob_man} amounts to a simple reweighting of the derivative terms in Equation  \ref{eq:sob_Rd}. The Sobolev spaces $W^{r,p}(M)$ for non-integer $r>0$ are, as in the case above, defined via interpolation theory. \cite{triebel2,fuselier} The non-integer spaces are crucial to the statement of trace theorems for Sobolev spaces, which are used in this paper to study the restrictions of functions that define certain RKH spaces. 
\vspace*{-.125in} 
\section{Reproducing Kernel Hilbert (RKH) Spaces}

\vspace*{-.1in} 

\label{sec:rkhs_review}
In this paper we make use of several properties of {\em real}, scalar-valued, RKH spaces.  Also, the analysis below is readily extended to real, vector-valued RKH spaces $\mathbb{H}:=H^k$ for $k\in \mathbb{N}$.  See \cite{bobade2} for the case where this is carried out in the context of consensus estimation. 

%
%
\vspace*{-.1in} 
\subsection{Basic Definitions and Constructions}
\vspace*{-.1in} 

An RKH space $H_X$ of functions that map a set $X \subseteq \mathbb{R}^d \rightarrow \mathbb{R}$  is defined in terms of a real-valued, continuous, symmetric, and positive type function  $\mathfrak{K}_X:X \times X\rightarrow \mathbb{R}$ that is referred to as the kernel underlying the RKH space. The subscript on $\mathfrak{K}_X$ is used to emphasize the set over which the kernel, as well as the functions in $H_X$ are defined.  When we say that $\mathfrak{K}_X$ is of positive type, this means that 
$
\sum_{i,j\leq 1}^N \alpha_i \mathfrak{K}_X(x_i,x_j) \alpha_j:=\alpha^T\mathbb{K}_{X,N}\alpha \geq 0
$
for all $\{ x_1,\ldots,x_N\} \subset \mathbb{R}^d$ 
and   $\alpha:=\{\alpha_1,\cdots,\alpha_N\}^T\in \mathbb{R}^N$, with the collocation matrix associated with $\{x_i\}_{1\leq i\leq N}$ defined as  $\mathbb{K}_{X,N}:=
[\mathfrak{K}_X(x_i,x_j)]\in \mathbb{R}^{N\times N}$. So all the collocation matrices of a  kernel of positive type are positive semidefinite. We say that the kernel is of strictly positive type if all of its collocation matrices $\mathbb{K}_{X,N}:=[\mathfrak{K}_X(x_i,x_j)]$ for distinct  points $\{x_k\}_{1\leq k\leq N}$  are strictly positive definite. 
The function $\mathfrak{K}_{X,x}:=\mathfrak{K}_X(x,\cdot)$ is known as the kernel function centered at $x\in X$, and a candidate for the  inner product of two such functions $\mathfrak{K}_{X,x},\mathfrak{K}_{X,y}$ is defined to be $(\mathfrak{K}_{X,x},\mathfrak{K}_{X,y})_{H_X}:=\mathfrak{K}_X(x,y)$ for all $x,y\in X$.  The RKH space $H_X$ is the closed finite span of the set of functions $\{\mathfrak{K}_{X,x}\ | \ x\in X\}$, that is, 
\vspace*{-.1in} 
\begin{align*}
    H_X:&=\overline{ \text{span}\{\mathfrak{K}_{X,x}\ | \ x\in X \}} \\
    &=  \left \{f:X \rightarrow \mathbb{R} \ \biggl | \ f=\lim_{N\rightarrow \infty} \sum_{i=1}^N \alpha_{N,i} \mathfrak{K}_{X, x_{N,i}}\right \},
\end{align*}

\vspace*{-.2in} 
where $\alpha_{N,i}\in \mathbb{R}$ and $x_{N,i}\in X $.
The closure above is taken with respect to the candidate inner product. The Hilbert space $(H_X,(\cdot,\cdot)_{H_X})$ above is also known as the native space induced by the kernel $\mathfrak{K}_X$.  It is well-known \cite{berlinet,saitoh,paulsen} that with this construction the reproducing property 
$
(\mathfrak{K}_{X,x},f)_{H_X}=f(x)
$
holds for all $f\in H_X$ and $x\in X$. Any Hilbert space $H$ is in fact a RKH space if all of the evaluation functionals that act on $H$ are in fact bounded operators from $H\rightarrow \mathbb{R}$. If it is further known  that if for some positive  constant ${\overline{\mathfrak{K}}}$ we have  $\sup_{x\in X}\mathfrak{K}_X(x,x)\leq \bar{{\mathfrak{K}}}_X<\infty$, then the evaluation operator $\mathcal{E}_{x}:H_X \rightarrow \mathbb{R}$ given by 
$
\mathcal{E}_{x}:=\mathcal{E}_{H_X,x}: f\mapsto f(x)
$ 
is a uniformly bounded linear operator since 
$
    |f(x)| = |\mathcal{E}_{x} f| = \left|(\mathfrak{K}_{X,x},f)_{H_X}\right|\leq  \sqrt{\mathfrak{K}_X(x,x)}\|f\|_{H_X}. 
$
This implies that $\|f\|_{C(X)}\lesssim \|f\|_{H_X}$, and therefore we have the  continuous inclusion $H_X \hookrightarrow C(X)$. We will only consider kernels ${\mathfrak{K}_X}$ on $X$ for which such a constant $\overline{{\mathfrak{K}}}_X$ exists. 
\vspace*{-.1in} 

Later in the paper we also make extensive use of the closed subspaces 
$H_\Omega:= \overline{\text{span}\{ \mathfrak{K}_{X,x} \ | \ x\in \Omega \} }\subseteq H_X$ for some subset $\Omega \subseteq X$. These spaces are important in understanding how the new PE condition are applied. One important fact is that we have the $H_X-$orthogonal decomposition $H_X=H_\Omega \oplus V_\Omega$
with $V_\Omega$ the kernel of the trace or restriction operator on the set $\Omega\subseteq X$,
$
V_\Omega:=\left \{ f\in H_X \ | \ R_\Omega f=\left. f \right|_{\Omega}=0\right \}.
$
That is, $f\in V_\Omega$ if and only if $f(x)=0$ for all $x\in \Omega$. This fundamental property follows from the analysis in \cite{saitoh} and \cite{berlinet}. Finally, in some cases when we specifically discuss  spaces derived from restrictions of functions $H_X$ to a subset $\Omega$, we denote these RKH spaces as $R_\Omega(H_X)$. 
%
%
\vspace*{-.1in} 
\subsection{Separation of Closed Sets by Reproducing Kernels}

\vspace*{-.1in} 

\label{sec:separation}
The current paper is interested in understanding how the use of a RKH space can make precise certain notions of convergence in adaptive estimation. We want to understand the {\em geometric} implications of the PE condition, that is, what it implies  about the trajectories of the dynamical system and the PE set. Essentially, we will ``test convergence" in $X$ of trajectories $x(t)\rightarrow \bar{x}$ by the condition that $f\bigl(x(t)\bigr)\rightarrow f(\bar{x})$ for all $f\in H_X$. As we will see, it can be important for understanding persistence that the space $H_X$ contain enough functions to separate, in a certain sense, the points of $X$.  Here an example can illustrate the the problem.  It is known that it is always possible to induce a metric $d_K$ associated with the kernel $\mathfrak{K}_X$ as described in \cite{devito}.
The problem is, our semiflows will be continuous with respect to some metric $d_X$, and the topology induced by $d_X$ may not be the same as that generated by $d_K$.  In fact it is easy to come up with kernels for which this is the case. As noted in Remark 1 of \cite{devito}, the bilinear kernel $k_X(x,y):=x^{\top}y$ for $x,y\in \mathbb{R}^d$ induces a RKH space $H_X$ for which the only subsets that can be separated are linear manifolds. 
In this specific case, $d_K$ induces a topology that is strictly coarser than the usual topology on $X:=\mathbb{R}^d$. Specifically, the metric $d_K$ can be used to discriminate convergence to a particular line through the origin, but not convergence to a point on that line. 

\vspace*{-.1in} 
We will see that some  useful geometric insights regarding the PE condition and positive orbits result if we do not allow the  kernel to induce such a coarse topology.  We would like the metric generated by the kernel to be equivalent with that on the state space. Reference \cite{devito} gives one example of a useful and simple  separation property.  An RKH space $H_X$ is said to separate a subset $A\subset X$ if for each $b\notin A$ there is a function $f\in H_X$ such that $f(a)=0$ for all $a\in A$ and $f(b)\neq 0$. This condition can be used to prove that $d_K$ and $d_X$ define the same topology. However, we will employ an even stronger condition, one that is well-suited to the construction of native spaces that contain  well-known  classes of  differentiable functions.   We assume  the existence of a rich family of bump functions in $H_X$. We say that $b_{r,x}:X\rightarrow [0,1]$ is a bump function on $(X,d_X)$ associated with the open ball $B_{r,x}=\{y\in X \vert d_X(x,y)<r\}$ provided that 1) $b_{r,x}=1$ on a neighborhood of $x$, and 2) $b_{r,x}$ is zero outside a compact set contained in $B_{r,x}$. It is immediate that if for any open set $B_{r,x}$, there is an associated bump function $b_{r,x}\in H_X$, then the RKH space $H_X$ separates the $d_X$-closed subsets of $X$. We say that the space $H_X$ contains a rich family of bump functions if it contains a bump function $b_{r,x}$ for each open ball $B_{r,x}$.  The construction of smooth bump functions on $X=\mathbb{R}^d$ is a classical exercise in analysis on manifolds,  \cite{lee}  pages 49--51. In practice, the RKH space $H_X$ (even when $X\not = \mathbb{R}^d$) will be selected so that it contains them. See the proofs below of Corollaries \ref{cor:cor1} and \ref{cor:cor2}. 

\vspace*{-.125in} 
\section{The RKH Embedding Method}

\vspace*{-.1in} 

\label{sec:adapest_review}
In this paper, we study a model problem of adaptive estimation for uncertain nonlinear systems governed by ordinary differential equations that have the form
\vspace*{-.1in} 
\begin{align}
    \dot{x}(t) = A x(t) + B f(x(t)), \qquad x(0) = x_0
\label{eq:eq_orig}
\end{align}

\vspace*{-.325in} 
with $A \in \mathbb{R}^{d \times d}$ a Hurwitz matrix, $B \in \mathbb{R}^{d \times 1}$, and $f: \mathbb{R}^d \to \mathbb{R}$. This equation is a special case of the general form in Equation \ref{eq:est_decomp}, with $g:\mathbb{R}^d\rightarrow \mathbb{R}^d:=Bf$. Methods for ensuring that this system of ODEs has local or global solutions are well-known, \cite{khalil}, and in this paper we always assume that for each $x_0\in X$ the equations  have classical solutions on $\mathbb{R}^+$. In this equation, it is assumed that the matrices $A$ and $B$ are known, but the (nonlinear) function $f$ is unknown. The adaptive estimation problem considered in this paper  uses the observations of the full state, $x(t)$ for all $t \geq 0$, to construct estimates $\hat{x}(t) \to x(t)$ and $\hat{f}(t) \to f$ as $t \to \infty$. While $f$ is unknown in our adaptive estimation problem, information about this function is reflected in the choice of an hypothesis space $H$ of functions to which $f$ belongs. Perhaps the most familiar choice of hypothesis space $H$ is one that is finite dimensional
$
    H_n :=  \{ f = \sum_{i=1,\ldots,n} \alpha_i \phi_i  \}
$
with $\{ \phi_i \}_{i=1}^n$ some fixed set of basis functions and $\phi_i: \mathbb{R}^d \to \mathbb{R}$ for $1 \leq i \leq n$. If we suppose for the moment that the unknown function $f = \sum_{i=1}^n \alpha^*_i \phi_i \in H_n$, then one canonical choice of an estimator is
\vspace*{-.1in} 
\begin{align*}
    \dot{\hat{x}}(t) &= A \hat{x}(t) + B \Phi^T(x(t)) \hat{\alpha}(t) \\
    \dot{\hat{\alpha}}(t) &= - \Gamma^{-1} \Phi(x(t)) B^T P(x(t) - \hat{x}(t))
\end{align*}

\vspace*{-.325in} 
with $P \in \mathbb{R}^{d \times d}$ the symmetric positive definite solution of Lyapunov's equation $PA + A^T P = -Q$ for a user-designed symmetric positive definite matrix $Q \in \mathbb{R}^{d \times d}$, and $\Gamma \in \mathbb{R}^{n \times n}$ symmetric and positive definite. When the errors in state $\tilde{x}(t) = x(t) - \hat{x}(t)$ and parameter errors  $\tilde{\alpha}(t) = \alpha^* - \hat{\alpha}(t)$ are defined, it can be shown directly that the errors satisfy the equations

\vspace*{-.275in} 
{\small 
\begin{align*}
    \begin{Bmatrix}
    \dot{\tilde{x}}(t) \\ \dot{\tilde{\alpha}}(t)
    \end{Bmatrix}
    =
    \begin{bmatrix}
    A & B \Phi^T(x(t)) \\
    -\Gamma^{-1} \Phi(x(t)) B^T P & 0
    \end{bmatrix}
    \begin{Bmatrix}
    \tilde{x}(t) \\ \tilde{\alpha}(t)
    \end{Bmatrix}
    +
   E(t) 
\end{align*}
}

\vspace*{-.275in}
 for $E(t):=\{Be(t),0\}^T$, with $e(t) = 0$ if $f\in H_n$. If it happens that $f \notin H_n$, then $e(t):= f(x(t)) - \sum_{i=1}^n \alpha^*_i \phi_i(x(t)) := f(x(t))-f_n(x(t))$ with $f_n$ a suitable finite dimensional approximation of $f$.  Precise conditions on the exponential stability of this system are a classical topic in adaptive estimation for uncertain ODEs. See \cite{naranna,narendra1977b} when $e(t)=0$. When $e(t)\not = 0$, see \cite{PoFar,IaSu,HC} for related  discussions of ultimate boundedness of errors. 
\vspace*{-.1in}

In this paper, we are interested in a class of dynamical systems where the unknown function $f$ belongs to the RKH space $H$. The  generic RKH space $H$ may be the full space $H_X$ or one of its closed subspaces $H_\Omega$ described in Section \ref{sec:rkhs_review}. The plant, estimator and the learning laws for this case can be expressed as
\vspace*{-.1in} 
\begin{align}
    \dot{x}(t) &= A x(t) + B \mathcal{E}_{x(t)} f, \\
    \dot{\hat{x}}(t) &= A \hat{x}(t) + B \mathcal{E}_{x(t)} \hat{f}(t), \\
    \dot{\hat{f}}(t) &= \Gamma^{-1} (B \mathcal{E}_{x(t)})^* P (x(t) - \hat{x}(t)),
    \label{eq:ode_Rd}
\end{align}

\vspace*{-.325in} 
where $x(t)$, $\hat{x}(t)$, $A$, $B$, and $P$ are defined as above. But the (nonlinear) functions $f$ and $\hat{f}(t)$ belong to the RKH space $H$ and $\mathcal{E}_{x}: H \to \mathbb{R}^d$ is the evaluation functional that is defined as $\mathcal{E}_{x} f = f(x)$ for all $x \in X$ and $f \in H$. Furthermore, the term $\Gamma \in \mathcal{L}(H,H)$ in the above equation is a self-adjoint, linear positive definite operator. The error equation analogous to the classical case shown above is given by
\vspace*{-.1in}
\begin{align}
    \begin{Bmatrix}
    \dot{\tilde{x}}(t) \\ \dot{\tilde{f}}(t)
    \end{Bmatrix}
    &=
    \begin{bmatrix}
    A & B \mathcal{E}_{x(t)} \\
    -\Gamma^{-1} (B \mathcal{E}_{x(t)})^* P & 0
    \end{bmatrix}
    \begin{Bmatrix}
    \tilde{x}(t) \\ \tilde{f}(t)
    \end{Bmatrix} \nonumber\\
    &= \mathbb{A}(t)
    \begin{Bmatrix} \tilde{x}(t) \\ \tilde{f}(t)\end{Bmatrix}.
    \label{eq:rkhs_err}
\end{align}

\vspace*{-.3in}
\noindent 
Note that the evolution of the above error equation is in  $\mathbb{R}^d \times H$ as opposed to on $\mathbb{R}^d \times \mathbb{R}^n$ in the classical adaptive estimator case. Some elementary conditions that guarantee the existence of solutions, as well as their continuous dependence on initial conditions, are given in \cite{bmpkf2017,bmpkf2017C}. In this paper we always assume that for each $x_0\in X$ the equations admit a classical  solution $t\mapsto (\tilde{x}(t),\tilde{f}(t))\in \mathbb{R}^d\times H$ for $t\in \mathbb{R}^+$. The following theorem, which simplifies considerably the analysis in \cite{bmpkf2017,bmpkf2017C}, shows that this is reasonable for many common choices of the RKH space $H$.
\vspace*{-.1in} 
\begin{thm}
Suppose that the RKH space $H$ is generated by a kernel $\mathfrak{K}_X:X\times X \rightarrow \mathbb{R}$ for which $\sup_{x\in X}\mathfrak{K}_X(x,x) \leq \overline{\mathfrak{K}}<\infty$. Then for each $(\tilde{x}_0,\tilde{f}_0) \in  \mathbb{R}^d \times H$ there is a unique solution of Equation \ref{eq:rkhs_err} in $C^1([0,\infty),\mathbb{R}^d\times H)$.
\end{thm}

\vspace*{-.225in} 
\begin{pf}
As discussed in Section \ref{sec:rkhs_review}, the hypotheses guarantee that the evaluation operator $\mathcal{E}_x:H\rightarrow H$ is linear and uniformly bounded in $x\in X$. It is immediate that 
$$
\lim_{t\rightarrow \infty} \frac{1}{t}\log^+(\|\mathbb{A}(t)\|)=0
$$
with $\log^+(\xi):=\max(0,\log(\xi))$. As discussed on page 211 of \cite{barreira} the governing equations have a unique global solution in time.
\end{pf}

\vspace*{-.225in} 
The definition of the PE condition proves sufficient for convergence of function estimates $\hat{f}(t)\rightarrow f$ generated by the RKH embedding method, much as in the conventional, finite dimensional case.  The analysis of convergence of parameters (ie, functions in our case) is notoriously long, so in this short paper we merely outline the proof in a special case. The full and lengthy details (for general $P$ and Hurwitz $A$) are given in \cite{jia_acc}. 

We say that a family of functions $\mathcal{F}$ over a set $S$ is uniformly equi-continuous if for each $\epsilon>0$, there is a $\delta_\epsilon>0$ such that for all $f\in \mathcal{F}$ and $a,b\in S$, $|a-b|<\delta_\epsilon \Rightarrow \|f(a)-f(b)\|<\epsilon$.
\vspace*{-.1in} 
\begin{thm}
Suppose that $P=I$, $A$ is negative definite, $\mathcal{F}=\{f(x(\cdot))|f\in H,\|f\|_H=1\}$ is uniformly equi-continuous, and the trajectory $\Gamma^+(x_0)$ is persistently exciting in the sense of Definition \ref{defn:PE_rkhs}. Then the solution $\tilde{x},\tilde{f}$ of Equations \ref{eq:rkhs_err} satisfy
    $\lim_{t\rightarrow \infty} \tilde{x}(t)=0$ and 
    $\lim_{t\rightarrow \infty} \tilde{f}(t)=0$.
\end{thm}

\vspace*{-.225in} 
\begin{pf}
The proof that $\tilde{x}(t)\rightarrow 0$ follows along lines that are entirely analogous to the classical or finite dimensional case, see \cite{bmpkf2017C} for the details when arguments are lifted to the infinite-dimensional state space $\mathbb{R}^d \times H$. The conclusion that $\tilde{f}(t)\rightarrow 0\in H$ follows immediately from Theorem 3.4 of \cite{bsdr1997}, provided that we can prove that there exists constants $T,\Delta,\delta,\gamma>0$ such that for each $t\geq T$ and $f\in H$ with $\|f\|_H=1$ there is an $s\in [t,t+\Delta]$ such that 
\begin{equation}
\|B\|_{\mathbb{R}^d}\left | \int_s^{s+\delta} E_{x(\tau)}fd\tau \right |=
\left \| \int_{s}^{s+\delta} BE_{x(\tau)} f d\tau \right \|_{\mathbb{R}^d} > \gamma.
\label{eq:eqaa}
\end{equation}
However, the condition above can be shown to be equivalent to the PE Definition \ref{defn:PE_rkhs} provided that the integrand is smooth enough to eliminate the possibility of certain ``rapid switching'' behavior. The equivalence of conditions as in Equation \ref{eq:eqaa} to those similar to  Definition \ref{defn:PE_rkhs} in the classical, finite dimensional case have been studied in great detail. See \cite{naranna87} for a detailed discussion with excellent illustrative examples of pathological rapid switching in the finite dimensional case. In the case at hand, Equation \ref{eq:eqaa} follows from the fact that $f(x(\cdot))\in \mathcal{F}$, a family of uniformly equi-continuous functions. The lengthy details of the proof can be found in \cite{jia_acc}. 
\end{pf}
\vspace*{-.125in} 
\section{Semiflows and Persistence of Excitation (PE)}

\vspace*{-.1in} 
\label{sec:PE_RKHS}
In this section we recall of few of the basic definitions of dynamical systems theory that will be essential to the analysis of this paper.  The aim is to be able to define persistence of excitation, not only for the model problem in Equations \ref{eq:est_decomp} or \ref{eq:eq_orig}, but for more  general evolutions on metric spaces. In particular we obtain a PE condition that can be applied to flows on Riemannian manifolds, which encompass  a few of our examples.  A continuous semiflow or semidynamical system on the complete metric space $(X,d_X)$  is  defined in terms of a continuous  semigroup $\{S(\tau)\}_{\tau \geq 0}$ on $X$ .  The manner in which systems of ODEs can generate such a semigroup, and thereby a semidynamical system is well-studied \cite{walker,saperstone}. The positive orbit $\Gamma^{+}(x_0)$ starting at $x_0$  defined to be the set 
\vspace*{-.1in} 
$$
\Gamma^+(x_0):=\bigcup_{\tau\geq 0} S(\tau)x_0 \subseteq X.
$$

\vspace*{-.325in} 
The positive limit set $\omega^+(x_0)$ associated with the initial condition $x_0$ is defined to be \vspace*{-.1in} 
$$
\omega^+(x_0):=\bigcap_{t\geq 0} \overline{\bigcup_{\tau\geq t}S(\tau)x_0},
$$

\vspace*{-.325in} 
which is equivalently expressed as 
\vspace*{-.1in} 
$$
\omega^+({x_0})=\left \{ 
y\in X \ | \ \exists t_k \rightarrow \infty \text{ s.t.}  \lim_{k\rightarrow \infty } S(t_k)x_0  \rightarrow y
\right \}.
$$
%
%
%
\vspace*{-.1in} 
\subsection{Persistence and Positive Limit Sets}

\vspace*{-.1in} 

The next few results illustrate simple and often  intuitive relationships between persistently excited sets, positive orbits $\Gamma^+(x_0)$, and the positive limit set $\omega^{+}(x_0)$.  We start with a simple result that illustrates an intuitive notion of what the new PE Definition \ref{defn:PE_rkhs} entails.
\vspace*{-.1in} 
\begin{thm}
\label{th:pe_pt}
Let $\mathfrak{K}:X\rightarrow \mathbb{R}$ be a monotone nonincreasing radial basis function and suppose  the associated kernel $\mathfrak{K}_X(x,y):=\mathfrak{K}(d_X(x,y))$ induces the RKH space $H_X\hookrightarrow C(X)$ for some fixed $\epsilon > 0$. Define the measurable sets 
$
I_{t,\epsilon}=\left \{ \tau \in [t,t+\Delta] \ \biggl | \ x(\tau)\in B_{\epsilon,x_\infty}
\right \}
$
for each $t\geq T_0$ with $B_{\epsilon, x_\infty}$ the open ball of radius $\epsilon$ centered at $x_\infty$. If the 
the Lebesgue measure $\mu$ satisfies
$
\mu(I_{t,\epsilon}) \geq 
\gamma_\epsilon  >0 
$
for some constant $\gamma_\epsilon$ for all $t\geq T_0$, then the singleton indexing set $\Omega:=\{x_\infty\}$ and the closed subspace $H_\Omega \subset H_X$ are  persistently excited in the sense of Definition \ref{defn:PE_rkhs}. In particular, if $x(t)\rightarrow x_\infty$, the set $\Omega:=\{x_\infty\}$ and closed subspace $H_\Omega$ are persistently excited. 
\end{thm}
\noindent Before proving the above theorem, let us unpack the above definition to understand the relatively straightforward underlying idea. The interval $I_{t,\epsilon}$ is the set of times contained in the interval $[t,t+\Delta]$ during which the trajectory $t\mapsto x(t)$ is within $\epsilon$ of the point $x_\infty$. This theorem says that if a trajectory  spends at least $\gamma_\epsilon$ amount of time in each interval $[t,t+\Delta]$  in the ball of radius $\epsilon$ centered at $x_\infty$, then  $\Omega=\{x_\infty\}$ and $H_\Omega$ are persistently excited.

\vspace*{-.225in} 
\begin{pf}
Without loss of generality, we assume that the kernel $\mathfrak{K}_X$ is normalized so that $\mathfrak{K}_X(x_\infty,x_\infty)=1$. By definition $H_\Omega:=\text{span}\{ \mathfrak{K}_{X,x_\infty}\}$ when $\Omega=\{x_\infty\}$, and for each $f\in H_\Omega$ with $f=\alpha \mathfrak{K}_{X,x_\infty}$ we have $\|f\|_{H_\Omega}=\alpha=f(x_\infty)$. 
Only the lower bound of the persistency definition is problematic, and we compute directly that 
\begin{align*}
    &\int_{t}^{t+\Delta} \left ( \mathcal{E}^*_{x(\tau)} \mathcal{E}_{x(\tau)} f,f \right )_{H_X} d\tau \\
    & \hspace{0.5in} = \alpha^2 \int_t^{t+\Delta} \mathfrak{K}_X^2(x(\tau),x_\infty) d\tau \\
& \hspace{0.5in} \geq \alpha^2 \left ( \min_{0\leq \xi \leq \epsilon} \mathfrak{K}^2(\xi) \right ) \mu(I_{t,\epsilon}) \geq
\gamma_\epsilon \mathfrak{K}_\epsilon \|f\|_{H_\Omega}^2
\end{align*}
with $\mathfrak{K}_\epsilon:=\min_{0\leq \xi \leq \epsilon} \mathfrak{K}^2(\xi)$. 
\end{pf}

\vspace*{-.225in} 
Theorem \ref{th:pe_pt} gives a direct interpretation of the persistency condition when we consider a singleton $\Omega:=\{x_\infty\}$ in terms of visitation to a neighborhood of $x_\infty$. It also suggests that there are many choices of kernels $\mathfrak{K}_X$ that induce PE spaces $H_\Omega$ for any convergent trajectory $x(t)\rightarrow x_\infty$. 
The monotonicity of the kernel in the above theorem is satisfied for a host of common choices of RKH spaces, see Chapter 9 of \cite{wendland} for the definition of completely monotone radial basis functions and kernels. This fact illustrates a significant difference with the conventional PE definition: there are many convergent trajectories that simply are not classically PE for a given set of regressors. We also note that if $X=\mathbb{R}^d$, there is a direct extension of this theorem for the finite set $\Omega:=\{x_{\infty,1},\ldots,x_{\infty,M}\}$, see Lemma 3.4  in \cite{kur95}. 

We begin our study of the geometric nature of PE sets by noting that the forward orbit is always dense in PE sets.
\vspace*{-.1in} 
\begin{thm}
\label{th:th1}
Let 
 $H_X$   be an  RKH space of functions    over the domain $X$ and suppose that this RKH space includes a rich family of bump functions. If the PE condition in Definition  \ref{defn:PE_rkhs} holds for a subset  $\Omega\subseteq X$, then the forward orbit $\Gamma^+(x_0)$ is dense in $\Omega$, $
 \Omega \subseteq \overline{\Gamma^+(x_0)}$. That is, we have 
$$
y \in \Omega \quad \Longrightarrow  \quad \exists \{t_k\}_{k\in \mathbb{N}} \ \text{   with   }\  \lim_{k\rightarrow \infty} S(t_k)x_0 \rightarrow y.
$$
\end{thm}

\vspace*{-.325in} 
\begin{pf}
Suppose to the contrary that there is an $y\in \Omega$ for which there is  no such convergent sequence.  This means that there is an  open ball $B_{r,y}$ such that $\Gamma^+(x_0)\bigcap B_{r,y}=\emptyset$. 
But since we have assumed there is a rich collection of bump functions, there is a bump function $b_{r,y}\in H_X$ that satisfies
$
b_{r,y}(y) = 1
$,
$
b_{r,y}(x) =0 \quad \forall x\not \in C_{r,y}
$
with $C_{r,y}\subset B_{r,y}$ a compact set. However, 
from Section \ref{sec:rkhs_review}, $H_X=H_\Omega \oplus V_\Omega$ with $V_\Omega=\{f\in H_X \ | 
\ f(x)=0,  \forall x\in \Omega \}.$ It follows that $b_{r,y}\in H_\Omega$. 
Since $\Gamma^+(x_0)\bigcap B_{r,x}=\emptyset$, the integral in Definition \ref{defn:PE_rkhs} is equal to zero
\vspace*{-.1in} 
\begin{align*}
    &\int_{t}^{t+\Delta} \left(\mathcal{E}^*_{x(\tau)} \mathcal{E}_{x(\tau)}b_{r,y},b_{r,y}\right)_{H_\Omega} d\tau \\
    & \hspace{1in} = \int_t^{t+\Delta}|b_{r,y}(x(\tau))|^2 d\tau=0
\end{align*}

\vspace*{-.325in} 
for each $t\geq T$. 
Since 
\vspace*{-.1in} 
$$\|b_{r,y}\|_{{H}_\Omega}:=\|P_\Omega b_{r,y}\|_{H_X}=\|b_{r,y}\|_{H_X}\gtrsim \|b_{r,y}\|_{C(X)}>0,
$$

\vspace*{-.225in} 
this is a contradiction of the PE property in Definition \ref{defn:PE_rkhs} and the theorem is proven.
\end{pf}

\vspace*{-.225in} 
Note that Theorem \ref{th:th1} does not require that the set of times $t_k\rightarrow \infty$.  Recall, on the other hand,  that the positive limit set $\omega^+(x_0)$  is contained in the closure of all accumulation points  of the orbit $\Gamma^+(x_0)$ for sequences of the form $\{S(\tau_k)x_0\}_{k\in \mathbb{N}}$, as $\tau_k\rightarrow \infty$. 
Next, we discuss a relationship of the positive limit set $\omega^+(x_0)$ and a PE space $H_\Omega$ over the indexing set $\Omega\subset X$ in Definition \ref{defn:PE_rkhs}. 

\vspace*{-.1in} 
\begin{thm}
\label{th:th2}
Let $H_{X}$ be the RKH space of functions over ${X}$  and  suppose that this RKH space  includes a rich family of bump functions. If the PE condition in Definition  \ref{defn:PE_rkhs} holds for $\Omega$, then 
$
\Omega \subseteq \omega^+(x_0).
$
\end{thm}

\vspace*{-.225in} 
\begin{pf}
The proof of this result is similar to the argument in Theorem \ref{th:th1}, so we only outline it. For an arbitrary $y\in \Omega$ we  build a  sequence $\{x(t_k)\}_{k\in \mathbb{N}}:= \{S(t_k)x_0\}_{k\in \mathbb{N}}$ such that 
$
\lim_{k\rightarrow \infty} t_k =\infty,
\lim_{k\rightarrow \infty} S(t_k)x_0 = y.
$
Pick the arbitrary $y\in \Omega$ and fix $r_0>0$.  Choose $t_0 \in (T,T+\Delta)$ such that $x(t_0)\in B_{r_0,y}$.  Such an $x(t_0)$ must exist. If such a time does not exist, we  could choose a bubble function $b_{r_0,y}$ on $B_{r_0,y}$ as in the last example such that $\|b_{r_0,y}\|_{H_\Omega}>0$,   for which  the integral $\int_{T}^{T+\Delta} |b_{r_0,y}(x(\tau))|^2 d\tau =0$ would follow from the condition that $B_{r_0,y} \cap \{x(\tau)\ | \tau\in (T,T+\Delta)\}=\emptyset$. This is a contradiction of the PE condition. We can then set $r_1=r_0/2$ and repeat this process seeking a $t_1\in (2T, 2T+\Delta)$ such that $x(t_1)\in B_{r_1,y}$, and so forth to generate $\{t_k\}_{k\in \mathbb{N}_0}$ with $t_k\rightarrow \infty$ and $\{x(t_k)\}_{k\in \mathbb{N}_0}$ with $x(t_k)\rightarrow y$. These sequences satisfy the desired conditions above, and we must have $y\in \omega^{+}(x_0)$. 
\end{pf}
%
%
%
\vspace*{-.1in} 
\subsection{Persistence of Excitation for Semiflows on Manifolds}

\vspace*{-.1in} 

A careful reading of the Definition \ref{defn:PE_rkhs} makes clear that it depends on the orbit $\Gamma^+(x_0)$ of a continuous semiflow $\{S(t)\}_{t\geq 0}$ on a complete metric space $(X,d_X)$, a subset $\Omega\subseteq X$, and an admissible kernel $\mathfrak{K}_X$ that defines the RKH space $H_X$ (and therefore also the closed subspace $H_\Omega$). Since it applies to subsets of complete metric spaces, it makes sense to consider much more general systems than the ODEs in the model Equations \ref{eq:est_decomp} or \ref{eq:eq_orig}. For instance, we have the following result for semiflows on manifolds, the case when the state space $X=\Omega=M$ in the PE Definition \ref{defn:PE_rkhs}. Note that below the semigroup $S(t)$ that defines the positive orbit $\Gamma^+(x_0)$ is defined on all of $M=X$. 
\vspace*{-.1in} 
\begin{cor}
\label{cor:cor1}
Suppose that $M$  is a  compact,   connected, $d$-dimensional  Riemannian manifold, and $\mathfrak{K}_M$ is kernel that induces a native space $H_M$ whose norm $\|\cdot\|_{H_M}$ is equivalent to that of the Sobolev space $W^{r,2}(M)$. If the orbit $\Gamma^+{(x_0)}$ persistently excites $H_M$, then $\omega^+(x_0)=M$. 
\end{cor}
\begin{pf}
\label{ex:ex1}
We first show that there are indeed such kernels $\mathfrak{K}_M$ that induce a native space $H_M \approx W^{r,2}(M)$. The Sobolev embedding theorem  on $\mathbb{R}^d$ states that $W^{r,2}(\mathbb{R}^d)$ is continuously embedded in $C(\mathbb{R}^d)$, $W^{r,2}(\mathbb{R}^d)\hookrightarrow C(\mathbb{R}^d)$ when $r>d/2$. As noted on page 1748 of \cite{HNW}, this fact can be used to conclude that $W^{r,2}(M)$ is continuously embedded in $C(M)$, $W^{r,2}(M)\hookrightarrow C(M)$. 
This means that we have  
$
|\mathcal{E}_{x}f|=|f(x)|\leq \|f\|_{C(\Omega)} \lesssim \|f\|_{W^{r,2}(M)}
$
for each $x\in M$ and $f\in W^{r,2}(M)$. In other words, 
each evaluation functional $\mathcal{E}_{x}: W^{r,2}(M)\rightarrow \mathbb{R}$ on $M$  is  bounded.  But  $W^{r,2}({M})$ is a Hilbert space; boundedness of all its evaluation functionals implies that  $W^{r,2}({M})$ is a RKH space. 
We define  the  Sobolev-Matern kernel $\mathfrak{K}_{M}^r$ of smoothness $r>0$ to be  the unique fundamental solution of the elliptic  differential  operator equation  $\sum_{1\leq \ell  \leq r}  (\nabla^{\ell})^*\nabla^\ell  \mathfrak{K}_{{M}}^r=\delta$ where $\nabla$ is the covariant derivative operator over the manifold ${M}$ and $\delta$ denotes the Dirac distribution. When we define the native space $H_M$ in terms of the Sobolev-Matern kernel $\mathfrak{K}_M^r$, we have  $H_M 
\approx W^{r,2}({M})$ for the chosen range $r>d/2$. The details of this analysis are given in   \cite{HNW} for the case when $M$ is a smooth Riemannian manifold that satisfies our standing assumptions on $M$,  or see reference  \cite{Opfer} for the special case $M:=\mathbb{R}^d$.    We next show  that the RKH space $H_M$  defined in this way contains a rich family of (smooth) cutoff or bubble functions. This proof is not surprising given what we know about Sobolev spaces on subsets of $\mathbb{R}^d$. One way to define $W^{r,2}(\Omega)$ is as the completion of $C^\infty(\Omega)$ in the Sobolev norm, so  the space $C^\infty(\Omega)$ is dense in $W^{r,2}(\Omega)$.  It is well-known that for any open ball contained  $\mathbb{R}^d$, there is a smooth cutoff function with compact support contained in that ball. This is a standard result in the study of manifolds and the construction of partitions of unity. \cite{lee}  It follows that the Sobolev space $W^{r,2}(\Omega)$ contains a rich family of bubble functions. The result extends more generally to Sobolev spaces $W^{r,2}(M)$ using the exponential map. The details of the proof are rather long, which we simply outline below.  (Particular examples of such a construction can be found in \cite{HNW} on page 1749 and again on page 1751 of the same reference.) If $\hat{f}$ is a cutoff function on a ball $B_{x,r}\subset \mathbb{R}^d$, it is possible to construct an associated cutoff function on the image $Exp_q(B_{x,r})\subset M$ under the  exponential map $Exp_q:T_qM\rightarrow M$ from $f:=\hat{f} \circ Exp_{q}^{-1}$.  Such an $f$ is always an element of $C^\infty(M)$ since $\hat{f}$ is just a smooth representation of $f$ with respect to a compatible $C^\infty$ chart. The only technical difficulty is showing that $f\in W^{r,2}(M)$. But this follows from Lemma 3.2 of \cite{HNW} which states that the exponential operator $Exp_q$ induces a map $g\rightarrow g\circ Exp_q$ that is  boundedly invertible from $W^{r,p}(Exp_q(\Omega))$ to $W^{r,p}(\Omega)$ for any measurable set $\Omega\subseteq \mathbb{R}^d$. Alternatively, we can argue that $W^{r,2}(M)$ is the completion of $C^\infty(M)$ (\cite{triebel2}, Section 7.4.5) with respect to the norm in Equation \ref{eq:sob_man}. We conclude that  $W^{r,2}(M)$ contains a rich family of (smooth) cutoff functions. If  the motion over the manifold $M$ satisfies the persistency condition in Definition \ref{defn:PE_rkhs}, then $M:=\omega^+(x_0)$. 
\end{pf}

\vspace*{-.1in} 
This example illustrates that the newly introduced persistency condition can be  applicable, {\em in principle}, to the study of certain evolutions over smooth   Riemannian manifolds. Still, the analysis in the example above is fairly abstract.  Perhaps more importantly, it is not a simple task to come up with a closed form expression for the Sobolev-Matern kernel. Of course this can be done for some standard manifolds like $\mathbb{R}^d$, the circle, or a torus, since the Sobolev-Matern kernels can be written down for these cases.  But it  is not readily accomplished  for some arbitrary  manifold ${M}$. The definition of the space $H_M\approx W^{r,2}(M)$ is intrinsic here: it depends on the (usually unknown) domain of the manifold  ${M}$, the atlas of charts used to define the manifold, and the covariant derivative operator  intrinsic to the manifold. 

\vspace*{-.1in} 
We next discuss how it is possible to come up with constructions of a kernel for ${M}$ that is extrinsic in the sense that it is defined by the restriction of some known kernel on a larger domain that contains $M$. This terminology is used in \cite{fuselier} that studies the approximation properties of spaces constructed in such a fashion. This line of attack is particularly  useful to the study of unknown or uncertain  dynamical systems via the RKH embedding method. The persistency of excitation condition is cast in terms of the kernel on the larger space in this case, which is assumed to have a known closed form expression. 
Carefully note that the forward orbit $\Gamma^+(x_0)$ in the following theorem is defined in terms of a semigroup $S(t):M\rightarrow M$, but $M$ is a proper subset of $X$. 
\vspace*{-.1in} 
\begin{cor}
\label{cor:cor2}
Let $M$ be an $m$-dimensional, smooth, compact, (regularly) embedded submanifold of $X:=\mathbb{R}^d$, and suppose that the $\{S(t)\}_{t\in \mathbb{R}^+}$ defines a $d_M$-continuous semiflow on $M$ with $d_M$ the metric on $M$. Denote by $\mathfrak{K}_X^r$ the Sobolev-Matern kernel on $X=\mathbb{R}^d$ for some $r>d/2$, define the kernel 
$
\mathfrak{K}_M(\cdot,\cdot)= \left. \mathfrak{K}^r_X \right|_M(\cdot,\cdot),
$
and denote by  $R_M(H_X)$ the RKH space generated by $\mathfrak{K}_M$. If the orbit $\Gamma^{+}(x_0)$ of the semiflow on $M$ persistently excites $R_M(H_X)$, then $M=\omega^+(x_0)$. 
\end{cor}

\vspace*{-.225in} 
\begin{pf}
 The Matern-Sobolev kernels over $\mathbb{R}^d$ are given for $r>d/2$ by
$
\mathfrak{K}_X(x,y)=\mathfrak{K}_X(\|x-y\|_{\mathbb{R}^d})
$
with
$
\mathfrak{K}(\xi)
=\frac{2^{1-(r-d/2)}}{\Gamma(r-d/2)}\xi^{r-d/2}\mathfrak{B}_{r-d/2}(\xi)
$
for all $x,y\in \mathbb{R}^d$ and $\xi:=\|x-y\|_{\mathbb{R}^d}$   with $\mathfrak{B}_{r-d/2}$ the Bessel function of order $r-d/2$. (\cite{fuselier}, page 1771 or \cite{HNRW}, page 1957) As in the last example, we have $H_X(\mathbb{R}^d)\approx W^{r,2}(\mathbb{R}^d)$ under the condition that $r>d/2$.
That the   candidate kernel  $\mathfrak{K}_{{M}}$ defined the restriction 
$
\mathfrak{K}_{M}(x,y):=\mathfrak{K} |_{M\times M}(x,y)
$
for all $x,y\in M$ is in fact an admissible kernel for a RKH space follows from standard results on RKH spaces,  \cite{berlinet} Section 4.2  and \cite{saitoh} Sections 2.2.1-2.2.2. 
At this point we do not yet have a rigorous notion of  exactly how smooth the restricted functions in $R_M(H_X)$ are, nor do we know whether the spaces $R_M(H_X)$ contain a rich set of cutoff functions.  
But from Lemma 4 of \cite{fuselier}, we know that  $R_M(H_X)=T(H_X)=T(W^{r,2}(\mathbb{R}^d)) $ where $T$ is the trace operator
$
T: f\rightarrow f|_{M}.
$
 From Proposition 2 of \cite{fuselier}, under the standing assumptions on ${M}$, the trace operator $T:f\rightarrow f |_{{M}}$ is a continuous operator from $W^{r,2}(\mathbb{R}^d)$ onto $W^{r-(d-m)/2,2}(M)$ for $r> (d-m)/2$ and $1\leq m\leq d$. In summary then, if we choose the kernel $\mathfrak{K}_{X}$ on $\mathbb{R}^d$ with a sufficiently large smoothness index $r$, we have 
 $
R_M(H_X)=
 T(W^{r,2}(\mathbb{R}^d))
 \approx W^{r-(d-m)/2,2}(M).
 $
 This set of equivalencies gives a precise notion of the smoothness or regularity of the restricted functions in the RKH space $R_M(H_X)$: the RKH space over $M$ is equivalent to the Sobolev space having smoothness $r-(d-m)/2$. 
 The remainder of the proof is now that same as in Corollary \ref{cor:cor1}.    
\end{pf}

\vspace*{-.225in} 
Note that the statement of persistence in Definition \ref{defn:PE_rkhs} is expressed in terms of the kernel $\mathfrak{K}_M:=
 \mathfrak{K}_X|_{M\times M}$, which can be used for computations since a closed form for $\mathfrak{K}_X$ is known.
 
 \section{Conclusions}
 This paper derives sufficient conditions for the convergence of function estimates in the RKH embedding method that are based on the recently introduced notion of persistently excited indexing sets $\Omega$ and subspaces $H_\Omega$ of an RKH space $H_X$. The paper establishes that persistently excited subsets are contained  as subsets of  the positive limit sets, if the RKH space has a rich collection of bump functions. We have also introduced both intrinsic and extrinsic methods for defining an appropriate RKH space in the event that the positive limit set is in fact certain types of smooth manifold. The extrinsic method seems particularly well-suited for the estimation of uncertain nonlinear systems since the form of the positive limit set is unknown. 
 
 The theoretical results of this paper establish that a reasonable choice of basis functions for practical finite dimensional approximations include radial basis functions (defined in terms of the kernel of the RKH space) that are centered on or near the positive limit set. It remains an open question as to how to devise versions of the RKH embedding strategy that adaptively selects the basis as estimation is carried out. 
 
 \vspace*{-.1in} 
 \begin{ack}
 \vspace*{-.1in} 
 Andrew J. Kurdila would like to acknowledge the support of the Army Research Office under the award \textbf{Distributed Consensus Learning for Geometric and Abstract Surfaces}, ARO Grant W911NF-13-1-0407.
 \end{ack}
\vspace*{-.1in} 
\bibliographystyle{unsrt}
\bibliography{perst}

\begin{thebibliography}{10}

\bibitem{kl2013}
Andrew Kurdila and Yu~Lei.
\newblock Adaptive control via embedding in reproducing kernel hilbert spaces.
\newblock In {\em 2013 American Control Conference}, pages 3384--3389. IEEE,
  2013.

\bibitem{bmpkf2017}
Parag Bobade, Suprotim Majumdar, Savio Pereira, Andrew~J Kurdila, and John~B
  Ferris.
\newblock Adaptive estimation for nonlinear systems using reproducing kernel
  hilbert spaces.
\newblock {\em Advances in Computational Mathematics}, 45(2):869--896, 2019.

\bibitem{bmpkf2017C}
Parag Bobade, Suprotim Majumdar, Savio Pereira, Andrew~J Kurdila, and John~B
  Ferris.
\newblock Adaptive estimation in reproducing kernel hilbert spaces.
\newblock In {\em 2017 American Control Conference (ACC)}, pages 5678--5683.
  IEEE, 2017.

\bibitem{hale_kocak}
Jack~K Hale and H{\"u}seyin Ko{\c{c}}ak.
\newblock {\em Dynamics and bifurcations}, volume~3.
\newblock Springer Science \& Business Media, 2012.

\bibitem{khalil}
Hassan~K Khalil.
\newblock Nonlinear systems.
\newblock {\em Upper Saddle River}, 2002.

\bibitem{walker}
JA~Walker.
\newblock Abstract dynamical systems and evolution equations.
\newblock In {\em Dynamical Systems and Evolution Equations}, pages 85--136.
  Springer, 1980.

\bibitem{narkud}
KS~Narendra and P~Kudva.
\newblock Stable adaptive schemes for identification and control.
\newblock {\em IEEE Trans. System. Man Cybernet, SMC-4}, 1974.

\bibitem{shimkin1987Persistency}
Nahum Shimkin and Arie Feuer.
\newblock Persistency of excitation in continuous-time systems.
\newblock {\em Systems \& control letters}, 9(3):225--233, 1987.

\bibitem{naranna87}
K.S. Narendra and A.M. Annaswamy.
\newblock Persistent excitation in adaptive systems.
\newblock {\em International Journal of Control}, 45(1):127--160, 1987.

\bibitem{Moore1992Functional}
JB~Moore, R~Horowitz, and W~Messner.
\newblock Functional persistence of excitation and observability for learning
  control systems.
\newblock {\em Journal of dynamic systems, measurement, and control},
  114(3):500--507, 1992.

\bibitem{Boyd1983On}
Stephen Boyd and Shankar Sastry.
\newblock On parameter convergence in adaptive control.
\newblock {\em Systems \& control letters}, 3(6):311--319, 1983.

\bibitem{sb2012}
Shankar Sastry and Marc Bodson.
\newblock {\em Adaptive control: stability, convergence and robustness}.
\newblock Courier Corporation, 2011.

\bibitem{naranna}
Kumpati~S Narendra and Anuradha~M Annaswamy.
\newblock {\em Stable adaptive systems}.
\newblock Courier Corporation, 2012.

\bibitem{IaSu}
Petros~A Ioannou and Jing Sun.
\newblock {\em Robust adaptive control}.
\newblock Courier Corporation, 2012.

\bibitem{PoFar}
Jay~A Farrell and Marios~M Polycarpou.
\newblock {\em Adaptive approximation based control: unifying neural, fuzzy and
  traditional adaptive approximation approaches}, volume~48.
\newblock John Wiley \& Sons, 2006.

\bibitem{bobade2}
Parag Bobade, Dimitra Panagou, and Andrew~J Kurdila.
\newblock Multi-agent adaptive estimation with consensus in reproducing kernel
  hilbert spaces.
\newblock In {\em 2019 18th European Control Conference (ECC)}, pages 572--577.
  IEEE, 2019.

\bibitem{HNW}
Thomas Hangelbroek, Francis~J Narcowich, and Joseph~D Ward.
\newblock Kernel approximation on manifolds i: bounding the lebesgue constant.
\newblock {\em SIAM Journal on Mathematical Analysis}, 42(4):1732--1760, 2010.

\bibitem{af2003}
R.~A. Adams and John Fournier.
\newblock {\em Sobolev spaces}, volume 140.
\newblock Elsevier, 2003.

\bibitem{lopez}
V{\'\i}ctor~Jim{\'e}nez L{\'o}pez, Gabriel~Soler L{\'o}pez, et~al.
\newblock Transitive flows on manifolds.
\newblock {\em Revista Matem{\'a}tica Iberoamericana}, 20(1):107--130, 2004.

\bibitem{triebel2}
Hans Triebel.
\newblock {\em Theory of Function Spaces, Volume 2}.
\newblock Birkhauser, 1992.

\bibitem{HNRW}
Thomas Hangelbroek, F~Narcowich, Christian Rieger, and J~Ward.
\newblock An inverse theorem for compact lipschitz regions in $\mathbb{R}^d$
  using localized kernel bases.
\newblock {\em Mathematics of Computation}, 87(312):1949--1989, 2018.

\bibitem{fuselier}
Edward Fuselier and Grady~B Wright.
\newblock Scattered data interpolation on embedded submanifolds with restricted
  positive definite kernels: Sobolev error estimates.
\newblock {\em SIAM Journal on Numerical Analysis}, 50(3):1753--1776, 2012.

\bibitem{berlinet}
Alain Berlinet and Christine Thomas-Agnan.
\newblock {\em Reproducing kernel Hilbert spaces in probability and
  statistics}.
\newblock Springer Science \& Business Media, 2011.

\bibitem{saitoh}
Saburou Saitoh and Yoshihiro Sawano.
\newblock {\em Theory of reproducing kernels and applications}.
\newblock Springer, 2016.

\bibitem{paulsen}
Vern~I Paulsen and Mrinal Raghupathi.
\newblock {\em An introduction to the theory of reproducing kernel Hilbert
  spaces}, volume 152.
\newblock Cambridge University Press, 2016.

\bibitem{devito}
Ernesto De~Vito, Lorenzo Rosasco, and Alessandro Toigo.
\newblock Learning sets with separating kernels.
\newblock {\em Applied and Computational Harmonic Analysis}, 37(2):185--217,
  2014.

\bibitem{lee}
John~M Lee.
\newblock {\em Introduction to smooth manifolds}.
\newblock Springer, 2001.

\bibitem{narendra1977b}
AP~Morgan and KS~Narendra.
\newblock On the stability of nonautonomous differential equations ̇x=a+b(t)x,
  with skew symmetric matrix b(t).
\newblock {\em SIAM Journal on Control and Optimization}, 15(1):163--176, 1977.

\bibitem{HC}
Naira Hovakimyan and Chengyu Cao.
\newblock {\em ℒ1 Adaptive Control Theory: Guaranteed Robustness with Fast
  Adaptation}.
\newblock SIAM, 2010.

\bibitem{barreira}
Luis Barreira and Claudia Valls.
\newblock Stability of nonautonomous differential equations in hilbert spaces.
\newblock {\em Journal of Differential Equations}, 217(1):204--248, 2005.

\bibitem{jia_acc}
Jia Guo, Sai~Tej Paruchuri, and Andrew~J Kurdila.
\newblock Persistence of excitation in continuously embedded reproducing kernel
  hilbert space.
\newblock In {\em (submitted to) 2020 American Control Conference (ACC)}. IEEE,
  2020.

\bibitem{bsdr1997}
J.~Baumeister, W.~Scondo, M.A. Demetriou, and I.G. Rosen.
\newblock On-line parameter estimation for infinite dimensional dynamical
  systems.
\newblock {\em SIAM Journal of Control and Optimisation}, 35(2):678--713, 1997.

\bibitem{saperstone}
Stephen~H Saperstone.
\newblock {\em Semidynamical systems in infinite dimensional spaces},
  volume~37.
\newblock Springer Science \& Business Media, 2012.

\bibitem{wendland}
Holger Wendland.
\newblock {\em Scattered data approximation}, volume~17.
\newblock Cambridge university press, 2004.

\bibitem{kur95}
AJ~Kurdila, Francis~J Narcowich, and Joseph~D Ward.
\newblock Persistency of excitation in identification using radial basis
  function approximants.
\newblock {\em SIAM journal on control and optimization}, 33(2):625--642, 1995.

\bibitem{Opfer}
Roland Opfer.
\newblock Multiscale kernels.
\newblock {\em Advances in computational mathematics}, 25(4):357--380, 2006.

\end{thebibliography}

\end{document}